\documentclass[12pt]{article}
\usepackage[all]{xy}
\usepackage{wrapfig}
\usepackage{hyperref}

\usepackage{epic}
\usepackage{amsmath}[1996/11/01]
\usepackage{amssymb,amsthm,amsfonts,latexsym,epsfig,graphics}
 \def\beql#1#2\eeql{\begin{equation}\label{#1}#2\end{equation}}

\advance \textheight by -1 \topmargin
\advance \textheight by -1 \headheight
\advance \textheight by -1 \headsep
\advance \textheight by -2 \footskip
\vsize = \textheight
\hsize = \textwidth
\newcommand{\operp}{\oplus}

\DeclareMathOperator{\Gal}{Gal}

\DeclareMathOperator{\Chi}{{\mathcal X}}

\DeclareMathOperator{\Irr}{Irr}

\DeclareMathOperator{\rad}{rad}

\DeclareMathOperator{\Char}{char}
\DeclareMathOperator{\disc}{disc}
\DeclareMathOperator{\GL}{GL}

\DeclareMathOperator{\Fix}{Fix}

\DeclareMathOperator{\Hom}{Hom}

\DeclareMathOperator{\id}{id}

\DeclareMathOperator{\SL}{SL}

\theoremstyle{plain}
\newtheorem{theorem}{Theorem}
\newtheorem{lemma}[theorem]{Lemma}
\newtheorem{assumption}[theorem]{Assumption}

\newtheorem{proposition}[theorem]{Proposition}
\newtheorem{corollary}[theorem]{Corollary}
\newtheorem{definition}[theorem]{Definition}
\newtheorem{remark}[theorem]{Remark}
\theoremstyle{remark}

\newtheorem{example}[theorem]{Example}

\numberwithin{theorem}{section}

\newcommand{\Z}{{\mathbb{Z}}}
\newcommand{\Q}{{\mathbb{Q}}}
\newcommand{\F}{{\mathbb{F}}}

\newcommand{\N}{{\mathbb{N}}}
\newcommand{\R}{{\mathbb{R}}}

\newcommand{\HH}{{\mathbf{H}}}
\newcommand{\C}{{\mathbb{C}}}
\newcommand{\trace}{\mbox{trace}}

\renewcommand{\em}{\sf}

\title{Orthogonal Stability
\footnote{This research is funded under 
Project-ID 286237555 -- TRR 195 -- by the
Deutsche Forschungsgemeinschaft (DFG, German Research Foundation).}}
\author{Gabriele Nebe and Richard Parker} 
\date{}
 
\begin{document}
\maketitle

 {\sc Abstract.} 
 A character (ordinary or modular) is called {\em orthogonally stable} 
 if all non-degenerate quadratic forms fixed by representations with 
 those constituents have the same determinant mod squares. 
 We show that this is the case provided there are no odd-degree orthogonal constituents. 
 We further show that if the reduction mod $p$ of an ordinary character is orthogonally stable, 
 this determinant is the reduction mod $p$ of the ordinary one. 
 In particular, if the  characteristic does not divide the group order, 
 we immediately see in which orthogonal group it lies. 
 We sketch methods for computing this determinant, and give some examples.
 \\
        MSC:  20C15; 20C20; 11E12; 11E57.
        \\
        {\sc keywords:} orthogonal characters of finite groups; decomposition matrices; blocks with cyclic defect groups.

\section{Introduction}
\label{sec:intro}

The ordinary and Brauer character tables of a finite group, 
with the (Frobenius-Schur) indicator given,
 specify the embeddings of the irreducible representations of the group 
in the classical groups over all finite fields,
 except that if the degree is even and the indicator is $+$,
 it does not specify in which of the two orthogonal groups the 
representation is embedded. 
 We show here that one further value - the {\em orthogonal discriminant} - 
added to the ordinary table enables this information 
also to be specified for all primes not dividing the group order.

A character determines the irreducible constituents, and we say that
a character is {\em orthogonally stable} if there is a representation with
that character fixing a non-degenerate quadratic form, and if none
of these constituents have odd degree and indicator $+$.  In this case,
the concept of orthogonal discriminant generalizes by taking proper
account of the characters with indicator $o$ and $-$.

The intended meaning of an orthogonally stable character is that it has a
 well-defined discriminant 
 (see Definition \ref{welldef}, Theorem \ref{discchi}, and Theorem \ref{OSeqOD}).
 Notice that, no matter which representation we take with this character 
over an arbitrary finite field, 
if it fixes a non-degenerate quadratic form at all, 
it fixes the one specified by the orthogonal discriminant.

For an ordinary character with Schur index 1, 
the orthogonal discriminant is related to the determinant (mod squares) 
of the fixed quadratic form, and \cite{DetChar}
 shows that it can be otherwise defined,
 still as an element of the character field mod squares, 
even if the Schur index is not $1$.  
Over a finite field, we denote the orthogonal discriminant by $O+$ or $O-$ depending on which orthogonal group is involved. 
 We show that the reduction mod $p$ of an ordinary character is $O+$
 (respectively $O-$) if the reduction is orthogonally stable and the ordinary 
orthogonal discriminant reduces to a square (respectively a non-square) 
 in the finite field (see Theorem \ref{unram} and Corollary \ref{primesnotdividingG}).

We then go on to consider ways of computing the orthogonal discriminant.  The ordinary one has {\it a priori} only a finite number of possible values - only primes dividing the group order and totally positive units mod squares can be involved -
 so {\it ad hoc} methods - restriction, tensors, reduction mod $p$ etc. - can often be used to compute them. 
In this note we focus on the use of decomposition matrices. 
The ordinary discriminant is not divisible by primes having an orthogonally 
stable reduction. If the character lies in a block of defect 1, then 
also the converse is true. Such arguments sometimes suffice to 
deduce all ordinary and most modular orthogonal discriminants directly
from the decomposition matrices, as illustrated for the 
first sporadic simple Janko group $J_1$. 
 There seems to be no guarantee, however, that these methods will always be sufficient, so we resort quite often to (at least some) explicit computations.

More may be true than we prove here.  We have yet to see an ordinary orthogonal discriminant divisible by an odd power of any prime dividing 2. 

We thank Thomas Breuer for
 his helpful comments improving the exposition of the paper 
and for checking consistency of the results 
of Section \ref{examples}.
 
\section{Overview of the paper} 

The paper starts with introducing 
 the relevant facts about quadratic
forms in Section \ref{qf}. 
The most important notion here is the discriminant of 
a quadratic form. We also recall  some rules how this discriminant
behaves with respect to orthogonal sums and extension and
restriction of scalars.
Section \ref{qfdvr}  compares to Section \ref{qffin} 
in the spirit of Brauer who uses group rings over discrete valuation
rings $R$ as a bridge between ordinary and modular representation theory. 
For each even dimension 
there are exactly two regular quadratic forms 
over $R$.
The reduction modulo the
maximal ideal is either of type $O+$ or $O-$, 
and the discriminant algebra either 
$R\oplus R$  or the unique unramified quadratic extension 
of $R$. This observation becomes crucial in Section \ref{Reduction}
where we compare orthogonal
discriminants of ordinary characters with the ones in their 
reductions modulo primes.

Section \ref{orthrep} gives the notion of orthogonal representations. 
These are $KG$-modules $V$ with a non-degenerate 
$G$-invariant quadratic form $Q$. The representation
hence takes values in the orthogonal group $O(V,Q)$. 
When the characteristic of $K$ divides the group order
the $KG$-module $V$ is 
not necessarily semisimple. Under mild assumptions in characteristic 2
we may replace it by 
some semisimple $KG$-module, the {\em anisotropic kernel},
having the same orthogonal discriminant
(see Section \ref{Secaniso}). 

The main notion of the paper, {\em orthogonal stability}, is defined in 
Section \ref{characters}. 
Orthogonally stable characters $\chi $ are exactly those that
have a well-defined orthogonal discriminant (Theorem \ref{OSeqOD})
in the sense of 
Definition \ref{welldef}. 
After briefly recalling the notions
of Brauer characters and reduction of ordinary characters,
we classify the orthogonally simple characters. 
These are orthogonal characters
not containing smaller orthogonal characters.  
They fall in three categories as given in Proposition \ref{orthsimple}.
Theorem \ref{dchisimple} expresses the orthogonal discriminants 
for two of the three categories (indicator $o$ and $-$)
in terms of character values. 
The main result of this paper is Theorem \ref{discchi} 
stating that orthogonally stable characters $\chi $ 
have a well-defined discriminant, $\disc (\chi )$, the 
{\em orthogonal discriminant of $\chi $}. 
To determine  $\disc (\chi )$ it 
suffices to compute the orthogonal discriminants
 of all the indicator $+$ constituents that occur in 
$\chi $ with odd multiplicity (see 
Proposition \ref{discchiOS} for an explicit formula).

The remaining two sections are devoted to showing that decomposition matrices
help to compute orthogonal discriminants for both modular and ordinary
characters. 
 Assume that $\Chi $ is an orthogonally stable ordinary character
 with character field $K$.
 If the reduction of $\Chi $ modulo a prime ideal $\wp $ in the 
 ring of algebraic integers of $K$ is 
 orthogonally stable, then this prime ideal does not divide 
 the orthogonal discriminant of $\Chi $
 (see Theorem \ref{unram} for a precise statement). 
 In particular all primes that divide 
 $\disc(\Chi )$ do also divide the group order leaving 
 only finitely many 
possible orthogonal discriminants (Corollary \ref{finOD}).
Subsection \ref{algo} sketches an algorithm to deduce
which of the square classes is the true orthogonal discriminant
by computing enough orthogonal discriminants of orthogonally stable
modular reductions. 
We illustrate these methods computing the orthogonal discriminants of 
the irreducible ordinary characters of the Held group $He$ from its
$p$-modular orthogonal discriminants and the decomposition matrices
for the primes $p$ dividing $|He|$.
For absolutely irreducible characters $\Chi $ of $p$-defect one,
the primes $\wp $ above $p$ do not divide the discriminant if and only 
if  the reduction of $\Chi $ modulo
$\wp $ is orthogonally stable (Theorem \ref{cycl}).
This observation is enough to compute the orthogonal discriminant 
of all ordinary characters for the group $J_1$ from its decomposition 
matrices that are available in GAP \cite{GAP}.

\section{Quadratic forms} \label{qf}

This section recalls some basic facts on quadratic forms.
For more details we refer the reader to the lecture notes \cite{Kneser}
or the more elaborate textbooks \cite{Scharlau}, \cite{OMeara}, or \cite{Tignol}.
The most important notion in this section is the discriminant of 
a quadratic form and its discriminant algebra as it is 
given in \cite[Section 10]{Kneser}.

Let $K$ be a field and $V$ a finite dimensional vector space over $K$.
A quadratic form is a map ${Q}:V\to K$ such that 
${Q}(ax) = a^2{Q}(x)$ for all $a\in K$, $x\in V$ and 
such that its {\em polarisation} 
\begin{equation}\label{polar} 
{B}: V\times V \to K, {B}(x,y) := {Q}(x+y) - {Q}(x)-{Q}(y) 
\end{equation} 
is a bilinear form. 
The determinant of ${Q}$ is the determinant of a 
Gram matrix of its polarisation. This is well defined modulo 
squares. 
 The
 quadratic form ${Q}$ 
 is called {\em non-degenerate}, if its polarisation is 
 non-degenerate, 
 i.e. $\det ({Q}) \neq 0$. 
 We then also call $(V,{Q})$ a non-degenerate quadratic space.

 \begin{definition} \label{discriminant}
	 If the characteristic of $K$ is not $2$ then 
	 the {\em discriminant} of $Q$ is 
	 $$\disc(Q):= (-1)^{{n}\choose{2}} \det(B) (K^{\times})^2 \in K/(K^{\times})^2 .$$
	 If $\delta \in K^{\times}$ represents $\disc(Q) $ then the
	 {\em discriminant algebra} of $Q$ is 
	 ${\mathcal D}(Q) := K[X]/(X^2-\delta )$.
 \end{definition}

 Note that the discriminant algebra of a non-degenerate quadratic form 
 is either a field extension of degree 2 of $K$ or isomorphic to $K\oplus K$.
 The latter is exactly the case if the discriminant of $Q$ is a square, 
 the {\em trivial discriminant}.

In characteristic 2 the bilinear form $B$ does not determine the quadratic
form $Q$ and one needs to replace the discriminant by the Arf invariant.
This is a class of $(K,+)/\wp (K) $, where
$$\wp (K) := \{ a^2+a \mid a\in K \} $$
is a subgroup of the additive group $(K,+)$. 
 
To simplify notation, we again call this invariant the 
discriminant of the quadratic form in characteristic 2. 

\begin{remark} \label{Char2dimeven}
	Let $K$ be a field of characteristic $2$ and $Q:V\to K$ a 
	non-degenerate quadratic form. 
Then the polarisation $B$ of $Q$ is symplectic and hence 
	$\dim(V) $ is even, say $\dim(V) = 2m$.
	\\ 
	By \cite[Definition 10.7]{Kneser} the quadratic space $(V,Q)$ 
	is the orthogonal sum of $2$-dimensional non-degenerate spaces
	$$(V,Q) = \perp _{i=1}^{m} \langle e_i,f_i \rangle $$
	with $B(e_i,f_i) =1 $. 
	Then 
	$$\disc(Q) := \sum_{i=1}^m Q(e_i) Q(f_i) + \wp(K) \in (K,+)/\wp(K) $$
	is well-defined. 
	We call $\disc(Q) $ the {\em discriminant} of $Q$ and 
	refer to $(K,+)/\wp(K) $ as the ``square classes'' of $K$.
	If $\disc(Q) = b+\wp (K)$ then we put  the
	{\em discriminant algebra} of $Q$ to be
	${\mathcal D}(Q) := K[X]/(X^2+X+b) $. 
	\\
	Note that $b\in \wp(K)$ if and only if the polynomial
	$X^2+X+b\in K[X]$ is reducible, if and only if 
	${\mathcal D}(Q) = K \oplus K$.
	In this case we say that $Q$ has {\em trivial discriminant}.
\end{remark}

There is a group structure on the set of quadratic $K$-algebras
(see \cite[Section 10]{Kneser}) 
so that one may treat discriminant algebras simultaneously in 
even and odd characteristic. 
However, we prefer to work with numbers, where one should keep in 
mind that multiplication of discriminants means addition of the 
representatives in characteristic 2.
 In this sense,
 for even dimensional quadratic spaces, the
discriminant is multiplicative with respect 
 to orthogonal direct sums:

 \begin{remark}\label{discmul}
	 If ${Q}$ and ${Q}'$ are two quadratic forms 
	 of even dimension then 
	 $$\disc ({Q} \perp {Q}') = \disc({Q}) \disc({{Q}'}) .$$
 \end{remark}

 \subsection{Restriction of scalars} \label{res} 

 Let $F$ be a field and $K$ be an extension field of finite degree, say $d$.
 Let $(V,Q)$ be an $n$-dimensional 
 quadratic space over $K$.
 Then restriction of scalars turns $V$ into an $nd$-dimensional space
 over $F$. If $T$ denotes the trace of $K/F$ then 
 $T\circ Q: V\to F$ is a quadratic form that is non-degenerate if 
 $Q$ is non-degenerate and $K/F$ is separable 
 (the latter condition is always fulfilled if $K$ is finite or of characteristic 0).
 Then  
 $$\det (T\circ Q) = N_{K/F}(\det (Q)) \disc(K/F) ^n $$ 
 where $\disc(K/F) $ is the field discriminant of $K$ over $F$.

 \begin{remark} \label{resdet} (cf. \cite[Lemma 2.2]{Milnor})
	 If $K/F$ is separable and 
 $n$ is even, then 
	 $$\disc (T\circ Q) = \left\{ \begin{array}{cc} 
		 N_{K/F} (\disc (Q)) & \Char(K) \neq 2 \\ 
		 T_{K/F}(\disc (Q)) & \Char (K) = 2 . 
	 \end{array} \right. $$
 \end{remark}

 \subsection{Hyperbolic forms} \label{hyp}

 For a finite dimensional $K$-vector space $V$  we put 
$$V^{\vee } := \Hom _K(V,K) = \{ f: V\to K \mid 
 f \mbox{ is } K-\mbox{linear } \} $$  to denote the dual space of $V$.
 Any non-degenerate bilinear form $B$ on $V$ yields an 
 isomorphism 
 \begin{equation}\label{Btilde} 
	 \tilde{B} : V\to V^{\vee}, v\mapsto (x\mapsto B(v,x)) 
	 .
 \end{equation}
 The {\em hyperbolic module} $\HH (V) $ is the quadratic $K$-space 
 $$(V \oplus V^{\vee } , {Q} ) \mbox{ with } {Q}((v,f)):=f(v) \mbox{ for all } v\in V, f\in V^{\vee } .$$
Then $\dim (\HH (V) ) = 2 \dim(V) $ is even and 
 the discriminant algebra ${\mathcal D}(\HH (V)) = K\oplus K$.
So hyperbolic modules always have trivial discriminant.

 \subsection{Quadratic forms over finite fields} \label{qffin}

 Let $\F_q$ denote the field with $q$ elements.
It is well known (see \cite[Theorem 11.4]{Taylor}) that 
there are two non-degenerate quadratic forms of dimension 2
over $\F_q$, the hyperbolic plane 
$\HH (\F_q)$ and the norm form
$N(\F_{q})$ of the quadratic extension of $\F_q$.
The underlying space for $N(\F_q)$ is $\F_{q^2}$, 
regarded as a 2-dimensional $\F_q$-space, and the quadratic form is
$${Q}(x):= x \overline{x} = x x^q = x^{q+1} .$$

On an $\F_q$-vector space
of even dimension $2m$ there are two isometry classes of non-degenerate
quadratic forms. 
These two forms can be distinguished by many 
invariants such as discriminants, Witt index, 
the number of isotropic vectors, and also the
order of the orthogonal group
 (see for instance \cite[Chapter IV]{Kneser}, \cite[Chapter 11]{Taylor}).
To fix notation we denote these two forms by 
$$Q_{2m}^+(q):=\HH (\F_q)^m \mbox{ and } 
Q_{2m}^-(q):=\HH (\F_q)^{m-1} \perp N(\F_{q}) .$$
The discriminant algebra of $Q_{2m}^+(q)$ is $\F_q \oplus \F_q $ 
whereas ${\mathcal D}(Q_{2m}^-(q)) = \F_{q^2}$.

\begin{definition} 
	We also use the symbols $O+$ and $O-$ to distinguish 
	between these two quadratic forms and sometimes write that 
	the discriminant of $Q_{2m}^+(q)$ is $O+$ and 
	$\disc( Q_{2m}^-(q) ) = O-$.
\end{definition}

In analogy to  Remark \ref{discmul} we get

\begin{remark} \label{discmulOD}
For $m,n\in \N $ we have 
	$Q_{2m}^{+}(q) \perp Q_{2n}^{+}(q) \cong Q_{2(n+m)}^+(q)$,  
	$Q_{2m}^{+}(q) \perp Q_{2n}^{-}(q) \cong Q_{2(n+m)}^-(q)$,  
	 and  $Q_{2m}^{-}(q) \perp Q_{2n}^{-}(q) \cong Q_{2(n+m)}^+(q)$.
\end{remark}

In the situation of Subsection \ref{res} the restriction of 
scalars from $\F_{q^d}$  to $\F_q$ does not change the type 
($O+$ or $O-$) of the quadratic space.

\begin{remark}
	Let $T:\F_{q^d} \to \F_q$ denote the trace.
	Then 
	$$T ( Q_{2m}^{+}(q^d) ) = Q_{2md}^+(q) \mbox{ and } 
	T ( Q_{2m}^{-}(q^d) ) = Q_{2md}^-(q) .$$
\end{remark}

\begin{proof}
	Let $Q$ be a non-degenerate even dimensional quadratic 
	form over $\F_{q^d}$ of discriminant $\delta $.
	\\
	First assume that $q$ is odd. Then
	by Remark \ref{resdet}  the discriminant of $T\circ Q $ is 
	$N_{\F_{q^d}/\F_q } (\delta ) $. 
	The norm is a group epimorphism between the 
	multiplicative groups of the two fields. 
	In particular it maps the unique subgroup $(\F_{q^d}^{\times })^2 $
	of index 2 in $\F_{q^d}^{\times }$ onto the 
	subgroup of squares in $\F_q^{\times } $.
	So $N_{\F_{q^d}/\F_q } (\delta ) $ is a square in $\F_q^{\times } $ if and only if $\delta $ is a 
	square in $\F_{q^d}^{\times }$
	\\
	Now assume that $q$ is even. 
	Then $\wp(\F_q) \leq (\F_q,+) $ is a subgroup of index 2
	and $T_{\F_{q^d}/\F_q} (\wp (\F_{q^d})) =\wp(\F_q) $.
	In particular the discriminant of $Q$ lies in $\wp (\F_{q^d})$
	if and only if the discriminant of $T\circ Q$ lies in 
	$\wp(\F_q)$.
\end{proof}

From the computations in characteristic 2 we may conclude the 
following observation, which is certainly well known:

\begin{corollary} 
	Let $q=2^d$ and $b\in \F_q$. Then $X^2+X+b \in \F_q[X]$ is irreducible,
	if and only if the trace $T_{\F_q/\F_2} (b) = 1$. 
\end{corollary}

\begin{remark}\label{ext+} (see \cite[Proposition 4.9 (a),(d)]{SinWillems})
Any $n$-dimensional 
quadratic space $(V,Q)$ over $\F_q$ extends to an $n$-dimensional quadratic 
space $(V\otimes \F_{q^d} , Q)$ over $\F_{q^d}$ by putting 
$Q(x\otimes a) := a^2 Q(x) $ for $x\in V$, $a\in \F_{q^d}$. 
For even degree extensions
this quadratic space over $\F_{q^d}$  is always of type $O+$. 
The type is, however, stable under odd degree field extensions.
	\begin{itemize}
		\item[(a)] If $d$ is even then
	$Q_{2m}^+(q) \otimes \F_{q^{d}}  \cong 
	Q_{2m}^-(q) \otimes \F_{q^{d}} \cong Q_{2m}^+(q^d) $.
		\item[(b)] If $d$ is odd then
	$Q_{2m}^+(q) \otimes \F_{q^{d}}  \cong Q_{2m}^+(q^d) $ 
	and 
	$Q_{2m}^-(q) \otimes \F_{q^{d}} \cong Q_{2m}^-(q^d) $.
	\end{itemize}
\end{remark}

\subsection{Quadratic forms over discrete valuation rings} \label{qfdvr}

In this section we recall some results on regular quadratic forms 
over local rings for the special situation needed in this paper. 
So let $R$ be a complete discrete valuation ring and $\pi $ a
generator of its unique maximal ideal. 
We assume that the residue field $F:= R/\pi R$ is a finite field with $q$ elements and 
that the field of fractions $K$ of $R$ has characteristic 0. 
A quadratic lattice $(M,Q)$ is a free $R$-module $M$ of finite rank 
together with an $R$-valued quadratic form $Q:M \to R$. 
The quadratic form is called {\em regular} if the map
$\tilde{B}$ from equation \eqref{Btilde} 
restricts to  an $R$-isomorphism between $M$ and its dual module $M^{\vee } = \Hom _R(M,R)$. This is equivalent to saying that the determinant of $Q$ is a 
unit in $R$. 

As shown in \cite[Satz 15.6]{Kneser} there is a bijection between 
isometry classes of regular quadratic $R$-lattices and 
non-degenerate quadratic forms  over the residue field $F$. 
In particular in our situation there are exactly two 
isometry classes of regular quadratic $R$-lattices $(M,Q)$ 
of even rank $2m$, corresponding to the two possible residue forms 
$(M/\pi M,\overline{Q}) \cong Q_{2m}^{+}(q) $ respectively $Q_{2m}^{-}(q)$.
We denote these regular quadratic lattices by $Q_{2m}^+(R)$ respectively
$Q_{2m}^{-}(R)$. Note that the norm form of the ring of integers 
in the unique unramified extension of degree $2$ of $K$ has as 
residue form $Q_2^-(q) = N(\F_q)$, so we get the following remark.

\begin{remark} \label{disc}
	The discriminant of $Q_{2m}^+(R) $ is $1$.
	If $\delta \in \disc (Q_{2m}^-(R)) $ then 
	$\delta $ is a unit in $R$ such that $K[\sqrt{\delta }]$ 
	is the unique unramified extension of degree $2$ of $K$. 
\end{remark} 

One important property of regular quadratic sublattices 
is that they split as orthogonal direct summands.

\begin{lemma} (see \cite[Satz 1.6]{Kneser}) \label{reg}
	Let $(M,Q)$ be a quadratic $R$-lattice and $N\leq M$ be some 
	$R$-submodule such that $(N,Q_{|N}) $ is regular. 
	Then $$(M,Q) = (N,Q_{|N} ) \operp (N^{\perp} , Q_{|N^{\perp }}) .$$
\end{lemma}

\subsection{Hermitian forms} 

Hermitian forms arise naturally when considering non self-dual characters 
as in Proposition \ref{orthsimple} (b). 
The formula from Proposition \ref{unitary} below is used in Theorem \ref{dchisimple} 
to find the orthogonal discriminant of such orthogonally simple 
characters.
For a Galois extension
$L/K$ of degree $[L:K]=2$  let $\overline{\phantom{x}} \in \Gal_K(L)$
denote the non-trivial Galois automorphism.
A Hermitian form ${H}$ on a non-zero $L$-vector space $V$
is a map
${H} : V\times V \to L$ that is $L$-linear in the
first argument and such that ${H}(y,x) = \overline{{H}(x,y)} $
for all $x,y \in V$. If $n:=\dim_L(V)$ then
restriction of scalars turns $V$ into a $K$-vector space of
dimension $2n$ and ${H}$ defines a
quadratic form ${Q}_H : V\to K, {Q}_H(x):={H}(x,x)$ for all
$x\in V$.

\begin{proposition} \label{unitary} (see \cite[page 350]{Scharlau})
	Let $(V,{H})$ be a non-degenerate Hermitian $L$-vector space. 
	\begin{itemize} 
		\item[(a)]
        Let $\Char (K) \neq 2$ and write $L=K[\sqrt{\delta }]$.
			Then $\disc ({Q}_H) = \delta ^n (K^{\times })^2$.
\item[(b)] If $K\cong \F_q$ with $q=2^d$, then 
	the discriminant of $Q_H$ is trivial if $n$ is even, and
	non-trivial if $n$ is odd. So
			${Q}_H \cong Q_{2n}^+(q)$ if $n$ is even, and 
			${Q}_H \cong Q_{2n}^-(q)$ if $n$ is odd.
	\end{itemize}
\end{proposition}


Note that \cite{Scharlau} also gives the Clifford invariant of
${Q}_H$ as well as the Arf invariant for general fields of characteristic 2.

\section{Orthogonal representations} \label{orthrep}

Let $G$ be a finite group and $K$ be a field. 
Any $KG$-module $V\neq \{ 0 \}$ defines a 
group homomorphism $\rho : G \to \GL(V)$.
Then $\rho $ is called a  $K$-representation of $G$.
Equivalence of representations is defined as isomorphism of 
$KG$-modules. 
The representation $\rho $ 
is called  {\em irreducible} if $V$ is a simple $KG$-module, 
i.e. $\{ 0 \}$ and $V$ are the only $G$-invariant submodules 
of $V$. 
The {\em trivial representation} is the map $\rho : G\to \GL(K)$, 
$\rho (g) = \id _K $ for all $g\in G$.

For any $KG$-module $V$ also the dual space $V^{\vee }$ is a 
$KG$-module. The corresponding representation $\rho ^{\vee }$ 
is called the dual representation of $\rho $. 
The representation $\rho $ is called {\em self-dual}, 
if $\rho ^{\vee} $ is equivalent to $\rho $. 

Given a representation $\rho $ or equivalently a $KG$-module $V$ we put
$${\mathcal Q}(V) = {\mathcal Q}(\rho ) := \{ {Q} : V \to K 
\mid Q \mbox{ is a quadratic form, } 
 {Q}(x\rho(g) ) = {Q}(x) \mbox{ for all } x\in V, g\in G \} $$
the space of $G$-invariant quadratic forms. 
We call $\rho $ an {\em orthogonal representation},
if ${\mathcal Q}(\rho )$ contains a non-degenerate quadratic form $Q$.
Then $(V,Q)$ is also called an {\em orthogonal} $KG$-module. 

\begin{remark} \label{firstorth}
\begin{itemize}
	\item[(a)] Let $(V,Q)$ be an orthogonal $KG$-module. 
Then the polarisation of  $Q$ 
		defines a $KG$-isomorphism $\tilde{B}$ (see equation \eqref{Btilde})
		between $V$ and its dual module $V^{\vee }$.
In particular orthogonal representations are self-dual.
\item[(b)]
	The set ${\mathcal Q}(V) $ is a vector space over $K$. 
	So if $Q\in {\mathcal Q}(V)$ then also $aQ\in {\mathcal Q}(V)$ 
	for all $a\in K$. 
	Clearly $\disc (aQ) = a^{\dim(V)} \disc(Q) $ 
	so for odd dimensional orthogonal representations
	the discriminants of the invariant quadratic forms 
	represent all  square classes of $K$ and hence odd 
		dimensional representations cannot have a
		well defined orthogonal discriminant.
\item[(c)] 
	If $V$ is a $KG$-module then the hyperbolic module 
		$\HH (V) = (V\oplus V^{\vee }, Q)$ from Section \ref{hyp}
		is an orthogonal $KG$-module. 
\end{itemize} 
\end{remark} 


\subsection{The anisotropic kernel}  \label{Secaniso}

If the characteristic of $K$ divides the group order, then not 
all representations are direct sums of irreducible ones. 
However, given an orthogonal $KG$-module $(V,Q)$ there is a $KG$-submodule $U\leq V$  with $Q(U) =\{ 0 \}$ such that the
sub-quotient $U^{\perp }/U$ is the orthogonal direct sum of simple 
orthogonal $KG$-modules, where here we need to assume that 
the trivial representation is not a constituent of $V$ if the characteristic
of $K$ is 2. 
The {\em anisotropic kernel} $U^{\perp}/U$ is uniquely determined up 
to $KG$-isometry and has the same discriminant as $(V,Q)$ (see Proposition \ref{aniso} below).

\begin{lemma}\label{triv}
	Let $\rho : G\to \GL(V)$ be an irreducible representation and
	 ${Q} \in {\mathcal Q}(\rho ) $  be a non-zero
	 $G$-invariant quadratic form on $V$. 
	Then either ${Q}$ is non-degenerate or $\Char(K) =2$,
	${Q}:V \to K $ is an $\F_2G$-module homomorphism and $\rho $ is 
	the trivial representation.
\end{lemma} 

\begin{proof}
	Assume that ${Q} $ is degenerate.
	Since its polarisation 
	${B}$ is $G$-invariant,
	the radical 
	$$\rad(B)  = V^{\perp} = \{ v\in V \mid B(v,w) = 0 \mbox{ for all } 
	w\in V \} $$
 is a $KG$-submodule of $V$.
	 As $V$ is simple and $\rad ({B} ) \neq 0$, we hence have 
	 $\rad(B) = V$ and hence ${B} = 0$. 
	Now ${B}(x,x)= 2{Q}(x)$ for all $x\in V$ so $Q\neq 0 = B $ 
	implies that $\Char(K)=2$. Moreover 
 the equation \eqref{polar} shows that 
	${Q}$ is a group homomorphism between the additive
	group of $V$ and $K$. The $G$-invariance of ${Q}$ implies
	that the kernel of this homomorphism is a submodule of the simple 
	$KG$-module $V$,
	so ${Q}$ is injective. Moreover for $x\in V$ and $g\in G$ 
	we have ${Q}(x+xg) = {Q}(x) + {Q}(xg) = 2{Q}(x) = 0$, so $x+xg \in \ker({Q}) = \{ 0\}$ and hence $x=xg$ for all 
	$x\in V, g\in G$. 
	This implies that $V$ is the trivial $KG$-module.
\end{proof}

\begin{definition}
	A submodule $U$ of the orthogonal $KG$-module $(V,Q)$ 
	is called {\em isotropic} if $Q(U) = \{0 \}$. 
	An orthogonal $KG$-module $(V,Q)$ is called
	{\em anisotropic} if it does not contain a non-zero 
	 isotropic $KG$-submodule.
\end{definition}

 Let $\rho :G\to \GL(V)$ be an orthogonal representation and
        ${Q} \in {\mathcal Q}(\rho )$ be a non-degenerate
        quadratic form.
        Let $U$ be a maximal
 isotropic $KG$-submodule of $V$.
Then $Q$ defines a quadratic form $\overline{Q}$ on $U^{\perp}/U$ by
$\overline{Q} (x+U ) := Q(x)$ for all $x\in U^{\perp }$.
The quadratic $KG$-module
        $(W,{Q}_0) := (U^{\perp}/U , \overline{Q}) $ is uniquely determined by
        $\rho $ and $Q$ up to $KG$-isometry and called the
	{\em anisotropic kernel} of $\rho $ (see for instance \cite[Lemma 4.2 (8)]{Dress}). 
Note that $V=\HH (U)$ is hyperbolic if and only if its anisotropic kernel
is $\{ 0 \}$. The discriminant of a 0-dimensional quadratic space 
is defined as 1. 

	As the referee pointed out the following proposition
	is also proved in \cite{DadeMon} and \cite{Thompson} who calls
	$(W,{Q}_0)$ the Witt kernel of $(V,Q)$.
	Both authors work with anisotropic bilinear 
	forms, for quadratic forms the extra assumption that $\rho $ 
	have no trivial constituent in characteristic 2 
	is needed (see \cite{witt2} for the general situation in characteristic 2).

\begin{proposition}\label{aniso}
        Assume that either
	 $\Char(K) \neq 2$ or the trivial
        representation is not a constituent of $\rho $. 
        Then the anisotropic kernel
        $(W,{Q}_0)$ of $(V,Q)$
        is the orthogonal sum of simple orthogonal $KG$-modules.
        Moreover $\disc (Q_0) = \disc(Q)$.
\end{proposition}

\begin{proof}
        If no simple submodule of $V$
        is isotropic, then put $W:=V$, $Q_0:=Q$.
        Otherwise there is a simple isotropic submodule
        $S\leq V$. Then ${Q}$ induces a quadratic form on $S^{\perp}/S$.
        Replacing $V$ by this smaller dimensional module we finally
        arrive at a module $(W,{Q}_0)$ for which the restriction of
        ${Q}_0$ to any simple submodule is non-zero.
        As we assumed that the trivial module is not a constituent
        of $V$ (and hence of $W$) when the characteristic of $K$ is $2$,
        Lemma \ref{triv} shows that the restriction of ${Q}_0$
        to any simple submodule $S$ of $W$ is non-degenerate.
        Hence $S$ splits as a direct orthogonal summand
        $W=S \operp S^{\perp }$. So $W$ is the orthogonal
        direct sum of simple quadratic $KG$-modules.
\end{proof}


\section{Characters} \label{characters}

We use capital letters $\Chi $ only for ordinary characters,
whereas $\chi $ may stand for either an ordinary or a Brauer character. 

\subsection{Brauer characters} \label{brauer}

In characteristic 0, the character of a $K$-representation $\rho : G\to \GL(V)$
is the map $\Chi_{\rho } : G \to K , g\mapsto \trace (\rho (g)) $ 
where the trace is the usual trace of an endomorphism.
One important property here is that two representations 
are equivalent if and only if they have the same character. 
\\
In positive characteristic, $p$, such a definition of character is not
information preserving. For instance the character of the 
sum of $p$ isomorphic modules would be identically 0. 
The definition of Brauer character can  be found in \cite[Chapter 2]{Navarro} and \cite[Chapter 15]{Isaacs}. 
Given a field $K$ of characteristic $p>0$ and a $K$-representation $\rho $ 
the Brauer character $\chi _{\rho }$ of $\rho $ is a map 
from the $p'$-conjugacy classes of $G$ to the field of complex numbers. 
For $g \in G$ of order $n$ not divisible by $p$ (or just $g\in G_{p'}$) 
the eigenvalues of $\rho (g)$ (in an algebraic closure $\hat{K}$ of $K$) are 
$n$-th roots of unity.
Conway polynomials provide 
one algorithmic method for agreeing on a consistent family of choices
 to identify
these roots of unity in $\hat{K}$ with complex roots of unity
(see  \cite[Introduction]{BrauerATLAS}). 
Then the Brauer character of $\rho $ assigns to $g$ the 
sum of these complex roots of unity.

The Brauer character of $\rho $ does in general not determine the representation
$\rho $ up to equivalence.
However, it determines the multiset of composition factors of $\rho $ up
to isomorphism.

\begin{definition}
	Let ${\Chi}$ be an ordinary  character of $G$.
	The {\em character field} $\Q ({\Chi}) $ is 
	the abelian number field 
	$\Q[{\Chi}(g) : g\in G] $ 
	generated over the rationals by all character values of ${\Chi}$.\\ 
	For a $p$-Brauer character $\chi $ of $G$ the {\em character field}
	$\F_p (\chi )$ is obtained by reducing the ring of integers of the  
	number field $L= \Q [\chi(g) \mid g\in G_{p'}] $ generated by the Brauer character values modulo any prime ideal that divides $p$.
	As $L$ is Galois over $\Q $, the finite field $\F_p (\chi )$ does not
	depend on the choice of the prime ideal. 
\end{definition}


\subsection{Reduction of characters} \label{reduction}

Let $K$ be a number field, $O_K$ its ring of integers, $\wp$ 
a non-zero prime ideal of $O_K$ with residue field $F:=O_K/\wp $. 
For a $K$-representation $\rho : G\to \GL(V)$ a {\em reduction mod} $\wp $ 
is obtained by choosing a $G$-invariant $O_K$-lattice $M$ in $V$,
i.e. a finitely generated $O_K$-submodule $M$ in $V$  
that contains a $K$-basis of $V$ and such that 
$M\rho(g) = M$ for all $g\in G$. 
Then $\overline{V}:=M/\wp M$ is an $FG$-module. 
The composition factors of $\overline{V}$ are independent of the 
choice of the $O_KG$-lattice $M$ in $V$  
(see for instance \cite[Theorem 15.6]{Isaacs}). 
The corresponding representation $\overline{\rho } : G\to \GL(\overline{V})$ 
is called a {\em reduction of $\rho $} modulo $\wp $. 
Whereas the representation $\overline{\rho }$ depends on the 
choice of the lattice $M$, its Brauer character only depends on $\rho $ and 
on $\wp $. 

\begin{remark}
Note that the character field of the Brauer character of $\overline{\rho}$
can be strictly contained in $F$. 
	The smallest example in \cite{ATLAS} are the two representations
	$\rho $ of degree $4$ of the group $2.S_5$, where the character field
	of $\rho $ is $\Q(\rho)= \Q[\sqrt{-3}]$. 
	The irrationality only occurs as the character value of the outer
	elements of order $6$, thus the reduction $\overline{\rho }$
	of $\rho $ modulo $2$ has 
	a rational Brauer character of degree $4$ and character
	field $\F_2(\overline{\rho }) = \F_2$. 
\end{remark}

Given an ordinary character $\Chi $ with character field $K$, 
the {\em reduction of $\Chi $} modulo $\wp $ is 
 the Brauer character of a 
reduction of a representation $\rho $ (over a suitable field extension) 
affording $\Chi $. 
We can compute the reduction of $\Chi $ modulo $\wp $ directly 
from the character table as the restriction of some Galois conjugate of 
$\Chi $ to the $p'$-classes.
For the computations of ordinary orthogonal discriminants from $p$-modular 
ones using the decomposition matrix it is crucial to carefully distinguish 
between the different $p$-modular reductions arising from the
different prime ideals that divide $p$ 
(see Section \ref{J1} for examples).

\begin{remark}
For any finite abelian extension $K$ of $\Q $
the Conway polynomials from Section \ref{brauer} identify a unique 
prime ideal $\wp _0$ of $O_K$ that contains the prime $p$ so that 
the reduction of $\Chi $ modulo $\wp_0$ is 
the restriction of $\Chi $ to the $p'$-classes of $G$. 
\\
For any prime ideal $\wp $ of $O_K$ that divides $p$, there is 
	$\sigma \in \Gal(K/\Q )$  such that $\wp = \sigma^{-1}(\wp _0)$.
	Then 
the reduction of $\Chi $ modulo $\wp $ is obtained as 
 the reduction of $\sigma (\Chi )$ modulo $\wp _0$. 
\end{remark}

\subsection{Orthogonally simple characters} 

In the following the notion ``character'' denotes either
an ordinary character or a Brauer character of the finite group $G$.

\begin{definition}
	A character $\chi $ is called {\em orthogonal}, 
	if there is an orthogonal representation with character $\chi $. 
	An orthogonal character $\chi $ is called 
	{\em orthogonally simple}, if $\chi $ is not 
	a proper sum of orthogonal characters.
\end{definition}

\begin{remark} 
	If $\chi $ is any character then the character
	$\chi ^{\vee }$ of the dual representation is just
	the complex conjugate of $\chi $. 
 Remark \ref{firstorth} (c) shows that 
	$\chi + \chi ^{\vee }$ is always an orthogonal character. 
\end{remark}

In \cite{BrauerATLAS} 
the well known Frobenius Schur indicator for
absolutely irreducible ordinary characters, 
is extended to cover absolutely irreducible Brauer characters. 
We combine the two definitions:

\begin{definition}\label{FrobeniusSchur}
	Let $\chi $ be an absolutely irreducible character 
	or Brauer character. 
	\begin{itemize}
		\item[(+)]
	If $\chi $ is orthogonal then the indicator of $\chi $ is $+$. 
\item[(+)]
	In characteristic $2$ also
	the trivial character is of indicator $+$. 
\item[(o)] 
	If $\chi $ is not self-dual then its indicator is $o$.
\item[(-)]
	In all other cases the indicator of $\chi $ is $-$.
	\end{itemize}
\end{definition}

In particular the absolutely irreducible characters of 
indicator $+$ are exactly those that are afforded by some representation
$\rho $ (over some extension field of the character field)
for which ${\mathcal Q}(\rho ) \neq \{0 \} $.

\begin{proposition} \label{orthsimple}
	Let $\chi $ be an orthogonally simple character.
	Then one of the following holds: 
	\begin{itemize}
		\item[(a)] $\chi $ is absolutely irreducible of indicator $+$. 
		\item[(b)] $\chi  = \psi + \psi ^{\vee }$ for some absolutely irreducible character $\psi $ of indicator $o$. 
		\item[(c)] $\chi = 2 \psi $ for some absolutely irreducible character $\psi $ of indicator $-$. 
	\end{itemize}
\end{proposition} 

\begin{proof}
	Let $\psi $ be an absolutely irreducible constituent of $\chi $. 
	As $\chi $ is self-dual, also $\psi ^{\vee }$ is a constituent of $\chi $. If $\psi \neq \psi ^{\vee }$ then $\psi + \psi ^{\vee }$ is an 
	orthogonal sub-character of $\chi $ and hence $\chi = \psi + \psi ^{\vee} $. If $\psi = \psi ^{\vee }$ and $\psi $ has indicator $+$ then 
	already $\psi $ is orthogonal and hence $\chi = \psi $. 
	If $\psi $ is self-dual but not orthogonal, then all orthogonal 
	characters containing $\psi $ also contain $2\psi $.
	As $2\psi = \psi + \psi^{\vee } $ is orthogonal we are 
	in situation (c).
\end{proof}

\begin{corollary}
	The only orthogonally simple characters of odd degree 
	are the absolutely irreducible indicator $+$ odd degree 
	characters that are not the trivial Brauer character in
	characteristic $2$.
\end{corollary}

\begin{definition} \label{welldef}
	Let $\chi $ be an orthogonal 
	ordinary or Brauer character and denote by $K$ 
	its character field. 
 We say that 
 ``the orthogonal discriminant of $\chi $ is well-defined'' if 
	there is a quadratic $K$-algebra ${\mathcal D}(\chi )  $ 
	with the following property:
	Given an orthogonal representation
	$\rho $ over some field $L$ and with 
	character $\chi $. 
	Then all non-degenerate forms  ${Q} \in {\mathcal Q}(\rho )$
	have discriminant algebra ${\mathcal D}({Q}) = L \otimes _K {\mathcal D}(\chi )  $. 
	\\
	The quadratic $K$-algebra ${\mathcal D}(\chi )  $ is then called 
 the {\em  discriminant algebra of $\chi $}.
 \\
	If ${\mathcal D}(\chi )  = K[X]/(X^2-\delta )$ (for $\Char (K) \neq 2$)
	respectively 
	${\mathcal D}(\chi ) = K[X]/(X^2+X+b)$ (in characteristic $2$) then 
	$$\disc(\chi ) = \left\{ \begin{array}{cc} \delta (K^{\times })^2 & \Char(K) \neq 2 \\
	b + \wp(K) & \Char (K) = 2 \end{array} \right. $$ 
		is called the {\em orthogonal discriminant} of $\chi $.
		\\
		Also for Brauer characters $\chi $ 
		we sometimes write $\disc (\chi ) = O+$ if 
		$\disc (\chi )$ is trivial  and 
		$\disc(\chi )=O-$ for the non-trivial square class.
\end{definition}

By Remark \ref{firstorth} (b) no character of odd degree 
has a well-defined discriminant.

\begin{theorem} \label{dchisimple}
	The orthogonal discriminant of an
orthogonally simple character 
of even degree is  well-defined.
\end{theorem} 

\begin{proof}
	Let $\chi $ be an orthogonally simple character
of even degree and denote by $K$ the character field of $\chi $.
We go through the three cases from Proposition \ref{orthsimple}: 
\begin{itemize}
\item[(a)] Here $\chi $ is absolutely irreducible 
of indicator + and even degree.  For number fields $K$ the short note \cite{DetChar} 
shows that the orthogonal discriminant of $\chi $ is well-defined. 
For finite fields $K$, there is an orthogonal representation $\rho $ over $K$ with character $\chi $. 
As $\chi $ is absolutely irreducible the dimension of the $K$-space ${\mathcal Q}(\rho )$ is 1, 
so ${\mathcal Q}(\rho ) = \{ a {Q} \mid a\in K \}$ for some non-degenerate ${Q}$. 
		Clearly this holds also if we extend scalars, ${\mathcal Q}(L\otimes \rho ) =
		\{ a Q \mid a\in L \}$. 
		So $\chi $ is orthogonally stable and $\disc(\chi ) = \disc ({Q}) $. 
\item[(b)] 
If $\chi $ is a $p$-Brauer character and $\F_p (\psi ) = \F_p (\chi )$, 
then ${Q}$ is hyperbolic and hence has trivial discriminant. 
Otherwise ${Q} $ comes from a Hermitian form and 
$\disc (\chi )$ can be read off from Proposition \ref{unitary}.
\item[(c)]
	For number fields $K$ \cite[Theorem B]{Turull} asserts that 
		$\disc(\chi ) = 1$ in this case. 
		If $K$ is a finite field, then $\psi $ is the character
		of a representation over $K$. If $\rho :G\to \GL(V)$ 
		is a representation
		affording the character $2\psi $ and ${Q}\in {\mathcal Q}(\rho )$ non-degenerate, then the restriction of the quadratic form 
		${Q} $ to any simple submodule of 
		$V$ is identically zero. In particular the form ${Q}$  is hyperbolic and hence 
		has trivial discriminant.
\end{itemize} 
\end{proof}

\begin{remark} \label{totallypositive}
	In characteristic 0 the orthogonally simple characters are 
	exactly the characters of the irreducible $\R$-representations. 
	In particular any orthogonally simple character is the 
	character of a representation $\rho $ over some  real 
	number field $L$. Moreover $B:=\sum_{g\in G} \rho(g) \rho(g)^{tr} $ 
	is an invariant symmetric bilinear 
	form that is {\em totally positive definite}, i.e. for all
	ring homomorphisms $\epsilon $
	of $L$ into the reals the form $\epsilon(B)$ is 
	positive definite. Hence  the orthogonal discriminant
	of an orthogonally simple character of degree $2m$ 
	is $(-1)^m d (K^{\times})^2 $ for some 
	totally positive $d $ in the character field $K$.
\end{remark}

\subsection{Orthogonally stable characters} 

	An orthogonally stable character is the 
	sum of even degree orthogonally simple characters:

\begin{definition}
	An orthogonal character 
	$\chi $ is called {\em orthogonally stable} 
	if all absolutely irreducible indicator $+$ constituents of 
	$\chi $ have even degree. 
\end{definition}

Our main result of this section is the following theorem:

\begin{theorem}\label{discchi}
Any orthogonally stable character 
has a well-defined orthogonal discriminant in the sense of Definition \ref{welldef}.
\end{theorem} 

The proof of this important theorem is split into two cases: 
Subsection \ref{Prooffinite} proves Theorem \ref{discchi} for finite fields 
and Subsection \ref{Proofodd} for arbitrary fields of characteristic not 2.

There are orthogonal representations $\rho $ such that 
all non-degenerate forms in ${\mathcal Q}(\rho )$ are isometric
without the character of $\rho $ being orthogonally stable, 
as the following example shows.

\begin{example}\label{exmod}
	Take $G$ to be the normaliser of a Sylow-3-subgroup of 
	$O_4^{+}(\F_3)$. 
	Then $G$ is the extension of $C_3\times C_3$ 
	by an elementary abelian $2$-group of order 8.
	The Brauer character $\chi $ of the natural $\F_3 G$-module 
	is the sum of 4 one-dimensional $\F_3$-characters with indicator $+$.
	In particular it is not orthogonally stable.
Explicit computations show that 
	the space of $N$-invariant quadratic forms is of dimension $2$ and 
all non-degenerate forms in this space are isometric to 
	$Q_{4}^+(3)$. 
\end{example}

However the notions of ``well-defined orthogonal discriminant'' and
``orthogonally stable'' are equivalent, as the next theorem shows.

\begin{theorem} \label{OSeqOD} 
	An orthogonal character  $\chi $ 
	has a well-defined orthogonal discriminant if and only if $\chi $ 
	is orthogonally stable.
\end{theorem}

\begin{proof}
	Theorem \ref{discchi} 
	show that an orthogonally stable 
	character has a well-defined discriminant. 
	To see the converse, assume that $\chi $ is not orthogonally stable. 
	Then there is an absolutely irreducible 
	indicator $+$ constituent $\psi $ of $\chi $ that has odd degree. 
	\\
	Assume first that  the underlying characteristic is not 2.
Let $L$ be some finite extension of the character field of $\chi $ such that 
	there are orthogonal 
	$L$-representations $\rho _1$ with character $\psi $ and
	$\rho _2$ with character $\chi -\psi $ and choose 
	non-degenerate quadratic forms 
	$Q_1\in {\mathcal Q}(\rho _1)$ and $Q_2\in {\mathcal Q}(\rho _2)$. 
	Then for any $a\in L^{\times} $ the form 
	$aQ_1 \perp Q_2 \in {\mathcal Q}( \rho _1 \oplus \rho _2 ) $ is a 
	non-degenerate quadratic form in the direct sum of the two representations. As the dimension of $\rho _1$ is odd the discriminants 
$$\{ \disc (aQ_1 \perp Q_2 ) \mid a \in L^{\times } \} = \{ a (L^{\times })^2  
	\mid a\in L^{\times } \}  $$ 
	yield all possible square classes in $L^{\times }$. 
	So there is no well-defined orthogonal discriminant of $\chi $.
	\\
	In characteristic 2, the only indicator $+$ odd degree character
	is the trivial character (see Remark \ref{Char2dimeven}). 
	So here
	 $\psi $ is the trivial character.
	By assumption $\chi $ is orthogonal, in particular it has even dimension, 
	so the trivial character occurs at least twice in $\chi $.
	Let $\F_q:=\F_2(\chi ) $ be the character field of $\chi $ and
	$(V,Q)$ be an orthogonal $\F_qG$-module with character $\chi-2\psi $.
	Then the orthogonal $\F_q G$-modules $Q_2^+(q) \perp (V,Q)$  
	and $Q_2^-(q) \perp (V,Q)$ 
	both afford the character $\chi $ 
	but have different orthogonal discriminants. 
\end{proof}

\subsection{Proof of Theorem \ref{discchi} for finite fields} \label{Prooffinite}

To prove Theorem \ref{discchi} for finite fields of characteristic $p$
we use the well known fact that 
representations can be realised over their character fields. 

\begin{lemma}\label{lemsimple}
	Let $V$ be 
	a simple $KG$-module with Brauer character $\chi $. 
	If $\chi $ is orthogonally stable 
	then all non-degenerate $G$-invariant quadratic forms on 
	$V$ have the same discriminant. 
\end{lemma}

\begin{proof}
	As $V$ is simple, the endomorphism ring of $V$ is a field, say $L$, 
	containing $K$, and $V$ is also an $LG$-module, say $V_L$, of dimension 
	$\dim(V) / [L:K]$. 
	Let $\psi $ denote the Brauer character of $V_L$ and 
	$\Gamma = \Gal (L/K)$ is the Galois group of $L$ over $K$.
	Then 
	$\chi = \sum _{\gamma \in \Gamma } \gamma (\psi )$ is a 
	sum of pairwise distinct absolutely irreducible Brauer characters.
	\\
	Let $Q:V\to K$ be a non-degenerate $G$-invariant quadratic form on $V$.
	\\
	If the indicator of $\psi $ is $+$, then $\psi (1)$ is 
	even, as $\chi $ is orthogonally stable. 
	Moreover there is a
	$G$-invariant quadratic form $Q' : V_L \to L$ such that 
	$Q= T_{L/K} (Q')$. In this case the discriminant of $Q$
	can be obtained from the well-defined discriminant of $\psi $ 
	using Remark \ref{resdet}. 
	\\
	If the indicator of $\psi $ is $o$, then there is a unique $\gamma _0\in \Gamma $ (of order 2) such that $\psi^{\vee } = \gamma _0(\psi ) $. 
	Let $F \leq L$ denote the fixed field of $\langle \gamma _0 \rangle $. 
	Then $\psi + \psi^{\vee} $ is an orthogonally simple $F$-character
	of the simple orthogonal $FG$-module $V_{F}$.
	Again $Q=T_{F/K} (Q')$ for some non-degenerate quadratic form
	on $V_F$. By Theorem \ref{dchisimple} all such forms $Q'$ have
	the same discriminant from which Remark \ref{resdet} yields the 
	discriminant of $Q$. 
	\\
	It is not possible that the indicator of $\psi $ is $-$ as 
	then $\psi $ would occur twice as a constituent of $\chi $
	contradicting the fact that $\chi $ is the character of a 
	simple $KG$-module.
\end{proof}

\begin{proof} (of Theorem \ref{discchi} for finite fields) 
	Let $\chi $ be an orthogonally stable  Brauer character
	and $K:=\F_p(\chi )$ denote its character field. 
	Let $L$ be an extension field of $K$ and 
	let $(V,Q)$ be an orthogonal $LG$-module with character $\chi $.
 As $\chi $ is orthogonally stable, the trivial character is not
        a constituent of $\chi $ and hence we can apply Proposition \ref{aniso}
        to find an anisotropic orthogonal $LG$-module  $(W,Q_0)$
        having the same discriminant as $(V,Q)$.
	Then $(W,Q_0)$ is the orthogonal sum of simple orthogonal
	$LG$-modules. By Lemma \ref{lemsimple} the orthogonal 
	discriminant of these simple orthogonal summands is well-defined
	by their Brauer character and so is the discriminant of $(W,Q_0)$
	(use Remark \ref{discmul})
	and hence the one of $(V,Q)$.
\end{proof}

From the proof we get the following procedure to compute the orthogonal
discriminant of an orthogonally stable Brauer character $\chi $. 

Put $K=\F_p(\chi )$. For an orthogonally simple constituent 
$\psi $ of $\chi $ put $F:= \F_p(\psi )$ and let $\Gamma_{\psi} = \Gal (F/(K\cap F))$.
Then $$\psi_K := \sum _{\gamma \in \Gamma_{\psi} } \gamma (\psi ) $$ 
is a sub-character of $\chi $.
   If $\disc(\psi ) = \delta (F^*)^2$ and $p\neq 2$, then
 $$\disc (\psi_K ) = (\prod _{\gamma \in \Gamma _{\psi}} \gamma (\delta ) ) (K^*)^2
        =: N_{K} (\disc (\psi )) .$$
        \\
        If $\disc(\psi ) = \delta + \wp(F)$ (so $p = 2$), then
	$$\disc (\psi_K) = (\sum _{\gamma \in \Gamma_{\psi} } \gamma (\delta ) ) + \wp(K)
        =: T_{K} (\disc (\psi )) .$$
	 With this notation we compute 

	 \begin{proposition}  \label{discchiOS}
	 If $\chi = \sum _{j=1}^n (\psi_j)_K$ then 
	 $$\disc (\chi ) = \left\{ \begin{array}{cc} 
			 \prod _{j=1}^n N_{K} (\disc (\psi _j))  & \Char(K) \neq 2 \\
			 \sum _{j=1}^n T_{K} (\disc (\psi _j))  & \Char(K) = 2 
		 \end{array} \right. $$
	 \end{proposition} 

\subsection{Proof of Theorem \ref{discchi} for characteristic $\neq 2$} \label{Proofodd}

For fields of characteristic not 2 we can use a different strategy to
compute orthogonal discriminants based on the notion of 
the discriminant of an involution (see \cite{DetChar} for a 
brief application to orthogonal discriminants of 
characters and \cite{Tignol} for an 
exhaustive treatment of discriminants of involutions). 
The proof uses the idea of \cite{DetChar} to deduce that 
the orthogonal discriminant lies in
 the character field $K$, even though the character
might not be the character of a  $K$-representation.

This proof also shows that one can compute the orthogonal discriminants 
of invariant forms intrinsically in the group algebra $FG$ over the 
prime field $F$, which is either $\F_p $ or $\Q $.

\begin{remark}  (see \cite[Lemma 2.1, Remark 2.3]{DetChar}) \label{E-}
	Let $K$ be a field and $B\in K^{n\times n}$ a symmetric non-degenerate
	matrix. 
	Then the adjoint involution of $B$ 
	$$ \iota _B (X) := B X^{tr} B^{-1}  \mbox{ for all } X\in K^{n\times n} .$$
	defines an $K$-algebra anti-automorphism of order $2$ on $K^{n\times n}$.
	The $(-1)$ eigenspace  of $\iota _B$ is
	$$E_{-}(B) := \{ X \in K^{n\times n} \mid \iota _B(X) = - X \} 
	= \{ BY  \mid Y = -Y^{tr } \} .$$
	In particular 
	\begin{itemize} 
           \item[(a)] There is $X\in E_-(B) $ with $\det(X) \neq 0$ 
		   if and only if $n$ is even.
	   \item[(b)] For any $X\in E_-(B)$ we have $\det(X) (K^{\times })^2 
		   = \det (B) (K^{\times })^2 $.
	   \item[(c)] A matrix $g\in \GL_n(K)$ is in the orthogonal group
		   $g \in O(B)$ if and only if $gBg^{tr} = B$, 
			so if and only if $\iota_B(g) = g^{-1}$.
			Then $g-g^{-1} \in E_-(B)$.
		\item[(d)] If $n=2m$ is even then 
			we define the {\em discriminant of $\iota _B$} 
			as $\disc(\iota _B) = \disc (B) = (-1)^m \det (X) (K^{\times })^2 $ 
			where $X$ is any invertible element of $E_-(B)$.
	\end{itemize} 
\end{remark}

We now consider the group algebra $FG$ of the finite group $G$ over
the prime field $F$. 
The group algebra comes with a natural involution 
$$\iota : FG \to FG , \iota  (\sum _{g\in G} a_g g ) = \sum _{g\in G} a_g g^{-1} .$$
Let 
$$\Sigma := \{ x\in FG \mid \iota (x) = -x \} = \{ \sum _{g\in G} a_g (g-g^{-1} ) \mid a_g \in F \} $$ 
denote the subspace of skew elements. 

If the characteristic of $F$ divides the 
group order, then $FG$ is not semisimple.  
So let $J(FG) $ denote the radical of this finite dimensional $F$-algebra.
Then 
$$FG/J(FG) = \bigoplus _{i=1}^h A_i \cong \bigoplus _{i=1}^h D_i^{n_i\times n_i } $$ 
is a semisimple $F$-algebra, i.e. the direct sum of matrix rings over division
algebras $D_i$. In finite characteristic all the $D_i=:K_i$ are fields. 
If $F=\Q $
we put $K_i:=Z(D_i) $ to be the center of $D_i$ and $\dim _{K_i}(D_i) =:m_i^2$.

Clearly the natural involution $\iota $ preserves the radical and hence
gives rise to an involution $\iota $ on $FG/J(FG)$ permuting the 
simple direct summands.

\begin{remark} \label{deltai}
Let $\chi _1,\ldots , \chi _h$  represent the Galois-orbits 
	(Frobenius-orbits for Brauer characters) 
 of the absolutely irreducible characters of $G$, suitably ordered, so 
that $\chi _i$ belongs to $A_i$. 
	Then 
$K_i=F(\chi_i)$ is the character field of $\chi _i$ and 
$\chi_i(1) = m_i n_i $.
\begin{itemize}
\item[(a)]
If the indicator of $\chi _i$  is $+$, then $\iota (A_i) = A_i$ 
and the restriction $\iota _i$ of
$\iota $ to $A_i$ is an orthogonal $K_i$-linear involution.
		If $\chi _i(1) $ is even, then by \cite[Corollary 2.8]{Tignol} 
		(see also \cite[Proposition 3.8]{DetChar})
		the algebra $A_i$ contains invertible elements 
		$\delta _i$ such that $\iota_i (\delta _i) = - \delta _i$.
		Then by Remark \ref{E-} 
		$$\disc (\iota _i) = (-1)^{\chi_i(1)/2} \det (\delta _i )
		(K_i^{\times })^2 .$$
		Note that for $D_i \neq K_i$ the determinant of $\delta _i$ 
		needs to be replaced by its reduced norm (see \cite[Section 9]{Reiner}).
\item[(b)] 
If the indicator of $\chi _i$ is $o$ then there is $i'$ such that 
		$\chi _i^{\vee} $ is Galois conjugate to $\chi _{i'}$.
		\\
		(b1) If $i=i'$ then $\iota (A_i) = A_i$ and the restriction 
		$\iota _i$ of
		$\iota $ to $A_i$ is a unitary involution. 
		Put $K_i^+:=\Fix _{K_i}(\iota ) $,  a subfield of index $2$ in 
		$K_i$. 
		Then there is some $\delta _i \in K_i=Z(A_i)$ such that 
		$K_i = K_i^+[\delta _i]$ and $\iota_i(\delta _i) = - \delta _i$ and again we have (see also Proposition \ref{unitary})
		$$\disc (\iota _i) = (-\delta _i^2)^{\chi_i(1)} ((K_i^+)^{\times} )^2 .$$
		(b2) If $i\neq i'$ (this never happens in characteristic 0) 
		then $\iota (A_i) = A_{i'} $ and $\iota $ yields a 
		hyperbolic involution on $A_i\oplus A_{i'}$. 
		Then $y=(\delta_i,\delta_{i'}) = (1,-1) \in A_i \oplus A_{i'}$ satisfies 
		$\iota(y) = - y$ and the discriminant of
		the restriction of $\iota $ to $A_i\oplus A_{i'}$ is
		a square in $K_i= K_{i'}$.
	\item[(c)] 
		If the indicator of $\chi _i$ is $-$, then 
		$\iota (A_i ) = A_i $ and the restriction of
$\iota $ to $A_i$ is a symplectic $K_i$-linear involution.
		Again by \cite[Corollary 2.8]{Tignol} there is 
		an invertible $\delta _i \in A_i$ such that 
		$\iota (\delta _i) = -\delta _i$.
		Note that $\chi _i$ occurs in any orthogonal 
		character with even multiplicity and hence 
		$\det(\delta _i)$ contributes to an even power 
		to the orthogonal discriminant.
	\end{itemize}
\end{remark}

We now reorder the summands of $FG/J(FG) $ such that 
$\chi _{s+1},\ldots , \chi _h$ have indicator $+$ and odd degree 
and such that all $\chi_1,\ldots ,\chi_s$ are either of even degree 
or not self-dual.

\begin{proposition}
	There is $\Delta \in \Sigma $ such that 
	$\Delta + J(FG) = (\delta _1,\ldots , \delta _s, 0,\ldots ,0)$.
\end{proposition}

\begin{proof}
	Choose $\delta \in FG$ such that 
	$\delta +J(FG) = (\delta _1,\ldots , \delta _s,0,\ldots , 0) $.
	Then 
	$$\iota (\delta ) + J(FG) = \iota (\delta + J(FG)) = - \delta + J(G) .$$
	Put $\Delta := \frac{1}{2} (\delta - \iota (\delta )) $.
	Then $\Delta \in \Sigma $ and $\Delta + J(FG) = \delta + J(FG) $.
\end{proof}

\begin{proof} (of Theorem \ref{discchi} for fields of characteristic $\neq 2$) 
	Let $\chi $ be an orthogonally stable character and $L$ 
	be some field containing the character field $F(\chi )$. 
	Let 
	$\rho $ be any $L$-representation with character $\chi $ and 
	$Q\in {\mathcal Q}(\rho )$ a non-degenerate invariant quadratic form. 
	Let $B$ be the polarization of $Q$. 
	For simplicity we choose matrices with respect to a basis
	that is adapted to a composition series of $\rho $. 
	Then all matrices in $\rho (LG) $ are block upper triangular matrices 
	where the blocks correspond to simple $LG$-modules. 
	As $\rho (g) B \rho(g)^{tr} = B$, the adjoint involution 
	$\iota _B$ maps $\rho (g) $ to $\rho (g^{-1})$. 
	In particular $\rho (\Delta ) \in E_-(B)$. 
	As $\chi $ is orthogonally stable, all diagonal blocks of $\rho (\Delta )$ are invertible and so is $\rho (\Delta )$. 
	Also the determinant of $\rho (\Delta )$ is the product of the
	determinants of the diagonal blocks of $\rho (\Delta )$ and hence
	uniquely determined by the composition factors of $\rho $. 
	By Remark \ref{E-}
	the determinant of $\rho (\Delta )$ is the determinant of $B$ 
	modulo squares. 
	In particular $\det (B)$ is independent of the chosen representation
	$\rho $ with character $\chi $. 
\end{proof}

\begin{remark}
	From the construction of  $\delta _i$  in Remark \ref{deltai} 
	one concludes that the formula in Proposition \ref{discchiOS} 
	can also be used to compute the discriminant of $\rho (\Delta )$.
\end{remark}

\section{Reduction of orthogonal representations}  \label{Reduction}

In this section we fix  a number field $K$  with
ring of integers $O_K$. 
For a prime ideal $\wp $ 
we denote by $K_{\wp }$ the completion of $K$ at $\wp $,
by $O_{\wp}$ its valuation ring and $\pi $ a generator of the
unique maximal ideal $\wp O_\wp = \pi O_{\wp }$.
Put $F:= O_{\wp}/\pi O_{\wp} \cong O_K/\wp $ to denote the residue field of $O_{\wp}$.

Let $\rho : G\to \GL(V)$ be an orthogonal $K$-representation of $G$ and
$Q\in {\mathcal Q}(\rho )$ be a non-degenerate $G$-invariant 
quadratic form on $V$ with polarisation $B$. 
Clearly $Q$ extends to a quadratic form on the completion 
$V_{\wp } = V \otimes _K K_{\wp }$. 
An $O_\wp $-lattice $M$ in $(V_{\wp},Q)$ is called {\em even} 
if $Q(M)\subseteq O_{\wp }$, and {\em integral} if $B(M,M) \subseteq O_{\wp}$, 
i.e. $M$ is contained in its  
{\em dual lattice}, $M^{\#} =\{ x\in V_{\wp } \mid B(x,M) \subseteq O_{\wp }\}$.
Clearly even lattices are integral and, as $2Q(v)=B(v,v)$, 
being integral and even is equivalent if $2$ is a unit in $O_{\wp}$.

\begin{assumption}\label{assumption}
Assume that either $2$ is a unit in $O_{\wp }$ or 
the trivial representation is not a constituent of the 
reduction of $\rho $ modulo $\wp $. 
\end{assumption}

\begin{lemma} \label{inteven}
	Under Assumption \ref{assumption} 
	all integral $O_{\wp} G$-lattices in $V_{\wp }$ are even.
\end{lemma} 

\begin{proof}
The statement is clear if $2$ is a unit in $O_{\wp }$.
So assume that $\Char (F) = 2$ and that $M$ is an
integral $O_\wp G$-lattice in $V_{\wp }$. 
Then $2Q(x) = B(x,x) \in O_{\wp }$ and 
$Q(x+y) + O_{\wp } = Q(x) + Q(y) + O_{\wp }$ for all 
$x,y \in M$. Hence 
	$$Q: M \to (\frac{1}{2} O_{\wp})/ O_{\wp } $$ 
is a $G$-invariant  homomorphism of abelian groups.
Its kernel $N$ is a $\rho(G)$-invariant $O_{\wp}$-sublattice 
	of $M$ and $G$ acts trivially on the $F$-module $M/(N+\pi M )$. 
As the trivial module is not a  constituent of the reduction 
	of $\rho $ modulo $\wp $, 
we have $N+\pi M = M$, so $N=M$ by the well known Nakayama Lemma. 
This shows that $Q(M) \subseteq O_{\wp} $ and hence $M$ is even.
\end{proof}

A {\em maximal even $O_{\wp } G$-lattice} is an even, 
$\rho (G)$-invariant lattice $M$ in $V_{\wp}$ such that 
no proper $\rho(G)$-invariant overlattice of $M$ is even. 

\begin{lemma} \label{maxlat}
	Assume Assumption \ref{assumption} and
	let $M$ be a maximal even $O_{\wp } G$-lattice.
	Then $\pi M^{\#} \subseteq M \subseteq M^{\#} $ and 
	$(M^{\#} , \pi  Q)$ is an even lattice. 
	Reduction modulo $\wp $ yields 
	two non-degenerate quadratic $FG$-modules 
	$(M^{\#} /M , \overline{\pi Q})$ and $(M/\pi M^{\#} , \overline{Q})$.
\end{lemma} 

\begin{proof}
	Clearly $M$ is integral, so $M\subseteq M^{\#} $. 
	If $\pi M^{\#} \not\subseteq M$, then we put 
	$X:= \frac{1}{\pi } M \cap \pi M^{\#} +M$. 
	This is clearly a proper $O_{\wp } G$-overlattice of $M$. 
	Elementary computations 
	show that $X$ is integral, and hence even (by Lemma \ref{inteven}), 
	contradicting the maximality of $M$. 
	Similarly $(M^{\#} , \pi Q)$ is integral and hence even.
	\\
	Moreover the radical of the reduction $\overline{Q}$ of $Q$ modulo $\wp $ on $M/\pi M$ is $\pi M^{\#} $, by definition. 
	Similarly one concludes that $(M^{\#} /M , \overline{\pi Q})$ 
	is non-degenerate.
\end{proof}

\subsection{Reduction of orthogonal characters}\label{reductionorth}

\begin{theorem} \label{unram}
	Let $\Chi $ be an ordinary orthogonal character, $K$ its character
	field, $\wp $ a prime ideal of $O_K$ and  $F=O_K/\wp $ 
	denote its residue field, a finite field of characteristic $p$.
	Assume that  the reduction $\chi $ of $\Chi $ modulo $\wp$ 
	is orthogonally stable. Then 
	\begin{itemize}
		\item[(a)] 
	$\Chi $ is orthogonally stable. 
\item[(b)]
	$\wp $ is not ramified in the discriminant algebra
			${\mathcal D} (\Chi )$.
 \item[(c)] 
	 If $\wp $ is inert in ${\mathcal D} (\Chi )$ 
	then $\chi $ has orthogonal discriminant $O-$.
\item[(d)] 
	Assume	that  $\dim _{\F_p (\chi)} (F) $ is odd.
			If $\Chi $ has trivial discriminant or 
			$\wp $ is split in ${\mathcal D} (\Chi )$ 
then $\chi $ has orthogonal discriminant $O+$.
	\end{itemize} 
\end{theorem}

The reader can think of a prime ideal $\wp $ being ramified in ${\mathcal D} (\Chi )$ 
as synonymous to the informal statement that $\wp $ divides the orthogonal 
discriminant (however if $2\in \wp $ then being unramified is a stronger 
condition, cf. Corollary \ref{mod2} for an example).
Also a prime ideal $\wp $ of $K$ is 
inert in the quadratic extension $L/K$ of number fields, if and only if the 
completion of $L$ at $\wp $ is the unique unramified degree 2 field extension of 
the $p$-adic number field $K_{\wp }$. 
If $\wp O_L = P_1 P_2 $, i.e. the prime ideal is split, 
then the completion of $L$ at $\wp $ is $L_{\wp } = K_{\wp } \oplus K_{\wp }$ 
and not a field. 
Note that all prime ideals are split in $K\oplus K$.

\begin{proof} 
	(a) is clear. 
	\\
	(b)
	If the $\wp $-adic Schur index of $\Chi $ is 1, then 
	there is a $K_{\wp }$-representation $\rho : G\to \GL(V_{\wp })$ 
	affording the
	character $\Chi $. 
	Otherwise there is a suitable ramified extension $L$ of $K_{\wp }$ 
	and an $L$-representation $\rho $ affording $\Chi $. 
	Let $O_{\wp }$ denote the ring of integers either in $K_{\wp}$ or
	in $L$ and $\pi $ a generator of the
	maximal ideal of $O_{\wp }$. Then $F=O_{\wp}/\pi O_{\wp }$.
	Choose a non-degenerate $Q\in {\mathcal Q}(\rho )$.

	The orthogonal stability of $\chi $ implies Assumption \ref{assumption}.
	Choose a maximal even $O_{\wp } G$-lattice $M$ in $V_{\wp}$.
	Then Lemma \ref{maxlat} gives a chain 
	$\pi M^{\#} \subseteq M \subseteq M^{\#} $ of even $O_{\wp }G$-lattices. 
	Now $M^{\#}/M$ and $M/\pi M^{\#}$ are orthogonal $FG$-representations 
	whose Brauer characters sum up to $\chi $.
	So both are orthogonally stable, in particular of 
	even dimensions. 
	
Take $v_1,\ldots , v_s \in M$ such that their images 
form a basis of $M/\pi M^{\#} $ and put 
$N:= \langle v_1,\ldots , v_s \rangle _{O_{\wp}} $.
Then $(N,Q_{|N})$ is a regular quadratic $O_{\wp}$-lattice of
discriminant, say, $\delta _1$. 
	By Lemma \ref{reg} we get 
 $M = N \operp N^{\perp }$.
 Lemma \ref{maxlat} yields that 
	 $(N^{\perp}, \frac{1}{\pi } Q _{|N^{\perp} } )$ is 
	regular.  If $\delta _2$ denotes its discriminant,
	then $\disc (Q) = \delta _1 \delta _2$.
	In particular 
	 Remark \ref{disc} says that 
	${\mathcal D} _{\wp} :=
	K_{\wp }[\sqrt{\delta _1 \delta _2}]$ is an unramified 
	extension of $K_{\wp}$.
	Moreover the orthogonal discriminant of the quadratic
	$F$-space 
$(N/\pi N,Q_{|N} )  \operp 
	 (N^{\perp}/\pi N^{\perp }, \frac{1}{\pi } Q _{|N^{\perp} } ) $ 
	 is $O+$ if ${\mathcal D} _{\wp}=K_{\wp}$ is of degree 1 and 
	 $O-$ if ${\mathcal D} _{\wp }$ is a quadratic extension of $K_{\wp}$.
	\\
	(c) Here ${\mathcal D}_{\wp }$ has degree 2 over $K_{\wp} $ and hence
	by Remark \ref{ext+}  $\dim _{\F_p (\chi )}(F) $ is 
	odd and the orthogonal discriminant of $\chi $ is $O-$. 
	\\
	(d) When ${\mathcal D} _{\wp }=K_{\wp} $ then 
	Remark \ref{ext+} (a) shows that it is only 
	possible to read off the orthogonal discriminant of $\chi $ 
	from the computations over $F$ if 
	 $\dim _{\F_p (\chi )}(F) $ is odd.
	 In this case the orthogonal discriminant of $\chi $ is $O+$. 
\end{proof}

In particular for dyadic primes, i.e. those that divide 2, Theorem \ref{unram} 
yields good restrictions on the discriminant. We highlight this 
for rational characters:

\begin{corollary} \label{mod2}
The discriminant of an orthogonally stable rational character $\Chi $
with an orthogonally stable reduction modulo $2$ is $1$ modulo $4$.
\end{corollary} 

\begin{proof}
	Let $d\in \Z$ be square free so that $d(\Q^{\times })^2=\disc(\Chi )$. 
	As $\Chi $ is orthogonally stable modulo $2$, 
	the prime $2$ is either split in $\Q[\sqrt{d}]$ (i.e. $d\equiv 1 \pmod{8}$) or inert in $\Q[\sqrt{d}]$ (i.e. $d\equiv -3 \pmod{8}$). 
\end{proof}

\subsection{Primes not dividing the group order} \label{ordinaryprimes}

The Brauer character table modulo a prime not dividing the
group order is exactly the ordinary character table. 
Given an absolutely irreducible 
ordinary character $\Chi $ of even degree and indicator $+$ 
and some prime $p$ not dividing the group order we aim to compute the 
orthogonal discriminant of the $p$-Brauer character $\chi = \Chi $. 
Let $K=\Q (\Chi )$ be the character field of $\Chi $, $O_K$ be 
its ring of integers. Then the Conway polynomials 
define a unique prime ideal $\wp _0$ of $O_K$ with 
$p\in \wp_0 $. Put $F:=O_K/\wp_0 $. 
Then $F$ is the character field of $\chi $, $\chi $ is
the reduction of $\Chi $ modulo $\wp_0 $, and 
$\chi $ is orthogonally stable. 
If $2\not \in \wp_0$ then Theorem \ref{unram} implies the following corollary:

\begin{corollary} \label{primesnotdividingG}
	\begin{itemize}
		\item[(a)] 
			There is  $d\in O_K\setminus \wp_0$ 
			such that  $d (K^{\times})^2 = \disc(\Chi )$.
\item[(b)]
$\disc (\chi ) = O+$ if and only if $d +\wp_0 \in O_K/\wp_0 $ is a square in $F$.
\item[(c)]
$\disc (\chi ) = O-$ if and only if $d +\wp_0 \in O_K/\wp_0 $ is a non-square in $F$.
	\end{itemize}
\end{corollary} 

As there are only finitely many square classes $\delta (K^{\times })^2$ 
of a given number field 
$K$ such that $K[\sqrt{\delta }] / K$ is unramified outside a 
finite set we hence conclude the following corollary.

\begin{corollary} \label{finOD}
	Let $P$ be a finite set of rational primes and $K$ be some number field.
	There is a finite set of square classes of $K$ 
	that contain the discriminants of all
orthogonally stable characters with character field $K$ 
	of groups whose orders are divisible only
by the primes in $P$.
\end{corollary}

\subsection{Computing the orthogonal discriminant
of an ordinary character}\label{algo}

Theorem \ref{unram} can be applied to compute the orthogonal 
discriminant of an orthogonally stable ordinary 
character $\Chi $
by computing enough orthogonally stable reductions.

To this aim, assume that we are given an orthogonally stable character 
$\Chi $ of some finite group $G$. 
Let $K:= \Q (\Chi )$ denote its character field and $O_K$ the ring of 
integers of $K$. 
Put 
$${\mathcal P} := \{ \wp \unlhd _{\max} O_K \mid \mbox{ 
reduction mod $\wp $ of } \Chi  \mbox{ is not orthogonally stable } \} .$$
Then ${\mathcal P}$ is contained in the set of prime ideals dividing 
$|G|O_K$. In particular it is a finite set. 

Put $\epsilon := (-1) ^{\Chi (1)/2} $, so that $\epsilon \disc(\Chi )$ 
is a totally positive square class in $K$ (cf. Remark \ref{totallypositive}). 
Compute 
 $\delta _1,\ldots , \delta _t \in O_K$ such that the products 
$\prod _{j\in I} \delta _j$, for $I\subseteq \{ 1,\ldots , t\}$ 
represent all
totally positive square classes
$\delta (K^{\times })^2$  for which $K_{\wp } [\sqrt{\delta }] $ is unramified for 
	all $\wp \not\in {\mathcal P}$.
We assume that $(\delta_i (K^{\times })^2 \mid i=1,\ldots , t )
$ 
is $\F_2$-linearly independent in $K^{\times }/(K^{\times })^2 $. 
For a non-dyadic prime $\wp \not\in {\mathcal P}$ we define
$a_{\wp} \in \F_2^t$ by $a_{\wp}(i) = 0$ if $\delta _i $ is a 
square in the residue field $O_K/\wp $ and $a_{\wp }(i) = 1$ otherwise. 
Choose non-dyadic primes 
$\{ \wp_1,\ldots, \wp_t \}$ such that the matrix 
$ A$ with columns $a_{\wp _i}$, so  
$$A=(a_{\wp_1} ,\ldots , a_{\wp _t}) 
\in \F_2^{t\times t} $$ has full rank $t$. 
Now compute the orthogonal discriminant $\epsilon d_i $ of $\Chi \pmod{\wp _i}$ for 
$i=1,\ldots , t $ and put $x_i := 1$ if $d_i=O-$ and $x_i:=0$ if $d_i = O+$. 
Put $y := x A^{-1}  $ and $I:= \{ i\in \{ 1,\ldots , t\} \mid y_i = 1 \}$. 
Then 
$$\delta := \epsilon \prod _{i\in I} \delta _i (K^{\times })^2 = \disc (\Chi ).$$

\begin{remark}
	In practice we also use all dyadic primes $\wp $ where $\Chi $ has
	an orthogonally stable reduction $\chi $ as these yield even more
	restrictions than the non-dyadic primes. 
Instead of working only with $\delta _1,\ldots, \delta _t$ we need 
	to test the ramification behaviour of $\wp $ in $K [\sqrt{\epsilon d}]/K$ 
	for all products $d$ of the $\delta _i$. 
	By Theorem \ref{unram} we can exclude those $d$ where $\wp $ is
	ramified and the orthogonal discriminant of $\chi $
allows us to further exclude those
	$d$ where $\wp $ is inert (if $\disc (\chi ) = O+$) 
	respectively split (if $\disc(\chi ) = O-$). 
\end{remark}

\begin{example} 
	The group $J_2$ has two absolutely irreducible 
	ordinary characters of degree $224$ 
	with character field $\Q[\sqrt{5}]$ that are exchanged under
	an outer automorphism. So $\disc(224a) = \disc(224b) =: d(\Q[\sqrt{5}]^{\times })^2$.
	The decomposition matrices in GAP \cite{GAP} tell us that 
	both characters are orthogonally stable modulo $2$ and $7$.
	Constructing the representations we obtain that 
	they yield embeddings into 
	$O_{224}^+(\F_4)$ respectively $O_{224}^+(\F_{49})$. 
	As $|J_2| = 2^73^35^27$ the only possibilities for $d$ are 
	$1,3,(5+\sqrt{5})/2$, and $3(5+\sqrt{5})/2$. 
	Now $1$ and $3$ are squares in $\F_{49}$, whereas $(5+\sqrt{5})/2$ 
	is not a square. As the reduction mod 7 has orthogonal discriminant +,
	we are left with $1$ or $3$ as possible discriminants. 
	The prime $2$ is ramified in the extension $\Q[\sqrt{3},\sqrt{5}]/\Q[\sqrt{5}]$, which excludes the possibility that $d=3$ just using 
	the fact that the reduction modulo 2 is orthogonally stable. 
	We hence conclude that the discriminants are squares,
	$\disc(224a) = \disc(224b) = 1$.
\end{example} 

\subsection{Cyclic defect}

In this section we assume that $\Chi $ is an even degree absolutely irreducible
orthogonal character in some $p$-block with cyclic defect group. 
We refer the reader to the original article \cite{DadeCyc} for details. 
By \cite[Theorem 1, Part 2]{DadeCyc}
the reduction $\chi $ of $\Chi $  is multiplicity free. 
In characteristic 2, the trivial character is the unique 
odd degree indicator + irreducible Brauer character (see Remark \ref{Char2dimeven}), so 
the fact that $\chi $ is multiplicity free and of 
even degree allows to conclude that $\chi $ is 
orthogonally stable and hence Theorem \ref{unram} applies. 
So the theory in this section is void for $p=2$
and hence we may and will assume that $p\neq 2$. 
In particular we may consider bilinear forms instead of quadratic forms, 
all integral lattices are automatically even, 
and a prime ideal $\wp $ is ramified in the discriminant 
field extension, if and only if there is a $G$-invariant 
maximal integral lattice $L$ such that $L^{\#}/L$ has odd dimension
over the residue field $O_{\wp}/\wp$.

For blocks of defect 1 the converse of Theorem \ref{unram} is true. 

\begin{theorem} \label{cycl}
	Let $p$ be an odd prime and let 
 $\Chi \in \Irr_{\C }(G)$ be an 
orthogonally stable character in a $p$-block with defect 1.
	Let $K$ denote the character field of $\Chi $ and
$\wp $ a prime ideal of $O_K$ that contains $p$. 
Then  the reduction  $\chi $ of $\Chi $ 
        modulo  $\wp $
        is  orthogonally stable if and only 
	if $\wp $ is unramified (i.e. split or inert) in ${\mathcal D}(\Chi)$.
\end{theorem}

\begin{proof}
	Let $K_{\wp }$ be the completion of $K$ at the prime $\wp $,
	$L_{\wp }$ be some unramified extension of $K_{\wp }$ with
	ring of integers $O_{\wp }$ and 
	$\rho : G\to \GL_n(L_{\wp })$ be some representation 
	affording the character $\Chi $. 
	Fix a non-degenerate quadratic form $Q\in {\mathcal Q}(\rho )$ 
	and choose $M \leq L_{\wp }^n$ to be some maximal even
	$\rho(G)$-invariant $O_{\wp }$-lattice.
	The assumption that the $p$-defect of $\Chi $ is one 
	implies that all $O_{\wp }G$-lattice in $L_{\wp}^n$ are
	linearly ordered by inclusion 
	(\cite[Theorem 11]{Brauer}). 
	As $M^{\#}/M$ and $M/\wp M^{\#} $ are self-dual modules, 
	their composition factors are ordered as 
	$(S_1,S_2, \ldots , S_s, S=S^{\vee} , S_s^{\vee} ,\ldots , S_2^{\vee } , S_1^{\vee} )$ for $M^{\#}/M$ and 
	$(T_1,T_2, \ldots , T_t, T=T^{\vee} , T_t^{\vee} ,\ldots , T_2^{\vee } , T_1^{\vee} )$ for $M/\wp M^{\#}$, 
	possibly without the self-dual composition factors $S$, $T$. 
	As the module $M^{\#}/\wp M^{\#}$ is multiplicity free, all 
	these composition factors are pairwise non-isomorphic. 
	The assumption that the reduction of $\Chi $ modulo $\wp $ is
	not orthogonally stable hence implies that 
	both self-dual composition factors $S$ and $T$ do occur and
	are of odd dimension. 
	In particular $M^{\#}/M$ has odd 
	dimension so the $\wp $-adic valuation of the discriminant $\delta $ 
	of $Q$ is odd. 
	As $p\neq 2$, this is equivalent to saying that $\wp $ ramifies
	in $K[\sqrt{\delta }] / K$. 
\end{proof}

\begin{example}\label{def2}
	Let $G=\SL_2(8)$ and $\Chi $ be the 
	rational absolutely irreducible character of degree 
	$\Chi (1) = 8$. 
	Then the reduction  of $\Chi $ modulo 3
	is the sum of two characters of degree 1 and 7, and hence 
	not orthogonally stable. 
	However, the discriminant of $\Chi $ is $1$. 
\end{example}

Example \ref{def2} is covered by the following general result 
for characters in blocks with cyclic defect groups. 

\begin{theorem}\label{cycl2}
	Assume that $\Chi $ is an indicator $+$ absolutely 
	irreducible ordinary character of even degree in a $p$-block of cyclic 
	defect. Denote by $K$ its character field,
	$O_K$ the ring of integers and $\wp $ a prime ideal of $O_K$ 
	containing the odd prime $p$.
	\begin{itemize}
		\item If $\Chi $ belongs to the exceptional vertex then 
			its reduction modulo $\wp $  is 
			not orthogonally stable, if and only if 
			$\wp $ is ramified in the discriminant 
			algebra ${\mathcal D} (\Chi )$.
		\item If $\Chi $ does  not belong to the exceptional 
			vertex then $\wp $ is ramified in the discriminant
			algebra ${\mathcal D} (\Chi )$ if 
			and only if the defect is odd and the reduction 
			of $\Chi $ mod $\wp $ is not orthogonally
			stable.
	\end{itemize}
	In particular for even defect  only the 
	characters belonging to the exceptional vertex may have 
	ramification at $\wp $ in the discriminant field extension.
\end{theorem}

\begin{proof}
	We keep the notation of Theorem \ref{cycl}. 
	If $\Chi $ belongs to the exceptional vertex, then by \cite[Theorem (VIII.3)]{Plesken}  the $G$-invariant $O_{\wp }$-lattices form a chain and
	we have the same situation as in the proof of Theorem \ref{cycl} 
	and hence the conclusion follows with the same argument. 

	Now assume that $\Chi $ does not belong to the exceptional vertex 
	and we also assume that its reduction mod $\wp $ is not 
	orthogonally stable. So there are exactly two constituents 
	$\chi _S$, $\chi _T$ of indicator $+$ and odd degree.
	The paper \cite{head} investigates 
	the radical idealizer process, that starts with the $O_{\wp }$-order
	$\Lambda $ spanned by the matrices in $\rho (G)$ and successively
	constructs involution invariant overorders until arriving at the
	hereditary order $H(\Lambda )$, the head order of $\Lambda $.
	As $H(\Lambda )$ is invariant 
	under the canonical involution of $G$, 
	its lattices hence form a chain that is invariant under taking duals.
	The explicit form of $H(\Lambda )$ given in 
	\cite[Theorem 3.15]{head} 
 allows to conclude that $\chi _S$ and 
	$\chi _T$ occur in the same constituent of $H(\Lambda )$ 
mod $\wp $ if and only if the defect is even. 
	The same argument as in the proof of Theorem \ref{cycl} 
	now applied to $H(\Lambda )$ 
	allows to conclude the statement.
\end{proof}

\section{Two examples} \label{examples}

\subsection{The first Janko group $J_1$} \label{J1} 

Let $J_1$ denote the first Janko group of order $2^3\cdot 3\cdot 5\cdot 7\cdot 11\cdot 19$. 
Then there are seven indicator $+$ absolutely irreducible ordinary 
characters of even degree: 
The two Galois conjugate characters $56a$, $56b$ with character field 
$\Q[\sqrt{5}]$, two rational characters of degree $76$ and 
three Galois conjugate characters $120a$, $120b$, $120c$ with 
character field $\Q[c_{19}]$, the unique subfield of degree 3
of the 19th cyclotomic field. 
All these three character fields $\Q, \Q[\sqrt{5}],$ and $\Q[c_{19}]$ 
have narrow class number 1, which means that any ideal of their ring of
integers has a totally positive generator that is unique up to 
multiplication by squares of units. 

The reductions modulo $2,3,5$ of all these seven characters stay absolutely 
irreducible, so the only primes that can divide the discriminants 
are those that divide $7\cdot 11 \cdot 19$. 
With GAP \cite{GAP} we compute the reduction modulo 7, 11, 19 as follows 
$$
\begin{array}{|r|r|r|r|r|r|r|r|} 
	\hline 
	\Chi & 56a & 56b & 76a & 76b & 120a & 120b & 120c \\ 
	\hline 
	\pmod{7} & 56a & 56b & 1+75 & 31+45 & 45+75 & 31+89 & 120 \\ 
	\pmod{11} & 56 & 7+49 & 27 +49 & 7+69 & 1+119 & 56 +64 & 14+106 \\
	\pmod{19} & 1+55 & 22+34 & 76a & 76b  & 43+77 & 43+77 & 43+77  \\
	\hline 
\end{array} $$

By Theorem \ref{cycl} the determinant of both characters of degree 
76 is $7\cdot 11 = 77$. 

For the non-rational characters the ``reduction modulo 7,11,19'' 
means the reduction modulo the prime ideal $\wp _0$ 
from Section \ref{reduction}. 
As in our cases, the residue fields are the prime fields,
we can identify the ideal $\wp _0$ using  the GAP-command 
$$x:=Int(FrobeniusCharacterValue(y,p));$$   where $y$ is a name of the 
irrationality (here $y=ER(5)$ respectively $y=EC(19)$) and $p$ the 
corresponding prime. 
We obtain
$$\begin{array}{|c|c|c|c|c|}
	\hline 
	y & \sqrt{5} & \sqrt{5} & c_{19}  & c_{19}  \\
	\hline 
 	p & 11 & 19 & 7 & 11  \\
x & 4 & 9 & 0 & 6 \\ \hline \end{array} $$ 
We get the factorization $11=p_{11} p'_{11}$ where 
$p_{11}=4-\sqrt{5}$ and $p'_{11}$ is its Galois conjugate $p'_{11}=4+\sqrt{5}$. 
Similarly
$19 = p_{19} p'_{19}$ with $p_{19}=(9-\sqrt{5})/2$, and hence 
$$\disc(56a) = p'_{11}p_{19} (\Q[\sqrt{5}]^{\times })^2 = (17-4\sqrt{5}) (\Q[\sqrt{5}]^{\times })^2 ,\ \disc(56b) = (17+4\sqrt{5}) (\Q[\sqrt{5}]^{\times })^2 .$$
For the field $\Q(c_{19})$ we have 
$c_{19} = z+z^8+z^{8^2}+z^{8^3}+z^{8^4}+z^{8^5} $ where $z=e^{2\pi i/19}$ is a primitive
19th root of unity.  
A generator $\sigma $ of the Galois group of $\Q(c_{19})$ acts on 
the cyclotomic field by squaring $z$ and hence maps $c_{19}$ to 
$c_{19}' := 4-c_{19}^2$. Note that $\sigma $ acts on the characters 
as the cyclic permutation $(120a,120c,120b)$. 
We compute 
$$\alpha =4+c_{19}-c_{19}', \ \beta = 20+10c_{19}-3c_{19}', \ 
\gamma =7+3c_{19}-c_{19}' $$
as totally positive elements of norm 7, 11, respectively 19. 
Replacing $c_{19}$ by $0$ maps $\alpha $ onto $0\in \F_{7}$ and 
similarly mapping $c_{19}$ to $6$ maps $\beta $ to $0\in \F_{11}$. 
So we get 
$$\disc (120a) = \alpha \sigma(\alpha) \beta \gamma (\Q(c_{19})^{\times })^2  = (29-18c_{19} -9c_{19}') 
(\Q(c_{19})^{\times })^2  $$
and $\disc (120c) = \sigma (\disc (120a)) $, 
$\disc (120b) = \sigma^2 (\disc (120a)) $. 

\begin{theorem}
	The discriminants of the ordinary absolutely irreducible 
	orthogonal characters of $J_1$ are 
$$
	\begin{array}{|@{}r@{}|@{}c@{}|@{}c@{}|@{}c@{}|@{}c@{}|@{}c@{}|@{}c@{}|@{}c@{}|} 
	\hline
	\Chi & 56a & 56b & 76a & 76b & 120a & 120b & 120c \\ 
	\disc(\Chi ) & 17-4\sqrt{5} & 17+4\sqrt{5} & 77 & 77 & 29-18c_{19}-9c_{19}' &  47+9c_{19} + 18c_{19}' &  38+9c_{19}-9c_{19}'   \\
	\hline
\end{array} $$
	where the characters and conjugacy classes are as in the atlas \cite{ATLAS}
	and modular atlas \cite{BrauerATLAS} respectively.
\end{theorem} 

\begin{remark} \label{notGalois} 
	The discriminant algebra 
	 ${\mathcal D}(\chi_{56a})$ 
	 is not Galois over $\Q$, its normal closure  
	$L={\mathcal D}(\chi_{56a})[\sqrt{\disc(56b)}]$ 
	has Galois group $D_8$, the dihedral
	group of order $8$. 
\\
	Similarly  $K:={\mathcal D}(\chi_{120a}) $
	 is not Galois over $\Q$. Its normal closure 
	$K[\sqrt{\disc(120b)}, \sqrt{\disc(120c)}]$ 
	has Galois group $C_2 \wr C_3$ of order $24$.
\end{remark}

\begin{remark}
	In
 \cite[Conjecture 3.9]{Craven} the author speculates that the discriminant
	field of any absolutely
        irreducible indicator $+$ character
        of degree $2\pmod{4}$ is an abelian number field.
        The sporadic simple O'Nan group of order
	$ 2^9\cdot 3^4\cdot 5\cdot 7^3\cdot 11\cdot 19\cdot 31$
	provides a counterexample
        to this conjecture:
        Let $\Chi $ be one of the
 two Galois conjugate absolutely irreducible characters of degree $169290$
        with character field $\Q[\sqrt{2}]$ and indicator $+$.
	Let $K:= \Q[\sqrt{2}][\sqrt{\disc(\Chi)}] $ denote the
	discriminant algebra of $\Chi $. 
	We claim that $K$ is not Galois over $\Q$. 
	With GAP we compute the decomposition matrix modulo $31 $
	which shows that 
 the reduction of $\Chi $ modulo one prime ideal, say $\wp _{31}$,
	dividing $31$
	is orthogonally stable. 
	The reduction modulo the other prime ideal, say $\wp' _{31}$, is 
	not orthogonally stable. So Theorem \ref{cycl}  tells us that 
	$\wp _{31}$ is not ramified in the discriminant field extension,
	whereas $\wp'_{31} $ is ramified. Therefore 
	the prime ideals dividing $31$ in $K$ 
	have different ramification behaviour so $K/\Q $ is not 
	Galois.
	\\
        Using orthogonal condensation methods we computed the
        discriminant of the two
	Galois conjugate characters as 
	$(-53\pm 36 \sqrt{2} )\Q [\sqrt{2}]^{\times }$.
\end{remark}

	\subsection{The Held group} 

The Held group $He$ is a sporadic simple group of order 
$2^{10}3^35^27^317$. 
With the meat-axe  \cite{meataxe} 
we construct all irreducible indicator $+$ 
even degree modular representations for the prime divisors of the group 
order. The tables below give those degrees followed by the orthogonal discriminants. 

{\bf Characteristic 2} 
$$
\begin{array}{|cc|cc|cc|cc|cc|}
	\hline
	246	& O+	&   246& 	O+	&   680& 	O-	& 1920	& O+	& 2008& 	O-\\
	4352& 	O+& 	4608& 	O+
	& 21504 &	O+ &   21504& 	O+ & & \\
	\hline
\end{array}
$$

{\bf Characteristic 3} 
$$
\begin{array}{|cc|cc|cc|cc|cc|}
	\hline
	1920&	O-&	6172&	O-	&6272&	O-	&7650&	O- &  14400	&O-\\
	\hline
\end{array}
$$

{\bf Characteristic 5}
$$
\begin{array}{|cc|cc|cc|cc|cc|cc|}
	\hline
	104&	O+	&  680&	O+	&1240&	O-	& 4080&	O+&	4116&	O-&	6528&	O+\\	7650 &	O+&
		9640	&O-&    10860&	O- &   11900&	O+ &   14400&	O- &   22050&	O+\\
	\hline
\end{array}
$$

{\bf Characteristic 7}
$$
\begin{array}{|cc|cc|cc|cc|cc|}
	\hline
	   50& 	O+	&   426& 	O+	&   798& 	O+	& 1072& 	O-& 	1700& 	O+\\ 	3654& 	O+& 	6154& 	O-
	& 6272& 	O-   &  13720& 	O-   &  23324& 	O+ \\
	\hline
\end{array}
$$

{\bf Characteristic 17}
$$
\begin{array}{|cc|cc|cc|cc|}
	\hline
	  680&	O+	&4080&	O+&	4352&	O-&	6528&	O+\\
	  7650&	O+&    10880&	O+&    11900&	O+
	&23324&	O+ \\
	\hline
\end{array}
$$

The following table gives the information on the
orthogonal discriminants in characteristic $0$. 
The first column gives the number of the ordinary 
absolutely irreducible character $\Chi $ 
followed by its degree $\Chi (1)$ and its discriminant $\disc(\Chi )$. 
For the prime divisors of the group order for 
which the reduction of $\Chi $ is orthogonally stable 
we also give the orthogonal discriminant $O+$ or $O-$ of this reduction
computed from the orthogonal discriminants in the tables before. 
For characters in blocks with defect 1 such that 
the reduction mod $p$ is not orthogonally stable we display $p$,
to indicate that $p$ is a prime divisor of the discriminant
(see Theorem \ref{cycl}). 
All character fields are rational, except for the two algebraic conjugate
characters of degree $21504$, where the character field is 
$\Q[\sqrt{21}]$. Here ``$.,O+$'' indicates that the reduction of 
the character number 30 modulo the prime ideal generated by $4-\sqrt{21}$ 
is orthogonally stable of orthogonal discriminant $O+$ and 
the one modulo $4+\sqrt{21}$ is not orthogonally stable and {\it vice versa}. 

$$
\begin{array}{|c|c|c|c|c|c|c|c|}
\hline
No & deg & disc & \pmod{2} & \pmod{3} & \pmod{5} & \pmod{7} & \pmod{17} \\
\hline
6 & 680  & 21 & O- & & O+ & & O+ \\
12 & 1920 & 17 & O+ & O- & O- & O- & 17 \\ 
13 &  4080 & 1 &  & & O+ & O+ & O+ \\ 
14 &  4352 & 105 & O+ & & &  & O- \\ 
15 & 6272 & 17 & O+ & O- & & O- & 17 \\
16 & 6528 & 1 & O+ & & O+ & & O+ \\
19 & 7650 & -1 & & O- & O+ & O- & O+ \\
22 & 10880 & 1 & O+ & & O+ & O+ & O+ \\ 
25 & 11900 & 21 & O- & & O+ & & O+ \\ 
26 & 13720 & 17 &  & & O- & O- & 17 \\ 
27 & 14400 & 17 & & O- & O- & & 17 \\
30 & 21504 & 357+68\sqrt{21} & O+ & & ``.,O+'' & & 17 \\
31 & 21504 & 357-68\sqrt{21} & O+ & & ``O+,.'' & & 17 \\
32 & 22050 & -119 & & O+ & O+ & 7 & 17 \\
33 & 23324 & 1 & & & & O+ & O+ \\
\hline
\end{array} $$ 
For the characters of numbers $\not\in \{ 30,31,33 \}$ the 
proofs are all similar, so let us give the one for number 14. 
Here the possible discriminants are in 
$\{1,3,5,7,15,21,35,105\}$. 
As the reduction mod 2 is orthogonally stable of discriminant $O+$, Corollary
\ref{mod2} tells us
 that the discriminant of $\Chi $ is 1 mod 8, leaving 1 and 105 as the
only possibilities. But the discriminant is not a square modulo 17, so 
$\disc (\Chi ) = 105$. 

Characters number 30 and 31 are algebraic conjugate and so are their
discriminants. We give the proof for number 30. 
Here the character field is $K=\Q[\sqrt{21}]$. 
Its ring of integers has class number 1 and 
a totally positive unit $u$ that is not a square, $u=(5+\sqrt{21})/2$. 
Note that $3u=((3+\sqrt{21})/2)^2$ is a square in $K$.
The prime ideals dividing $3 \cdot 5 \cdot 7 \cdot 17 $ are 
generated by 
$p_3 := (3+\sqrt{21})/2 $ of norm $-3$, 
$p_5 := 4+\sqrt{21} $ and $p_5' := 4-\sqrt{21} $ of norm $-5$, 
$p_7 = (7+\sqrt{21})/2 $ of norm $7$ and $17$ of norm $17^2$. 
By Theorem \ref{cycl} we know that 17 divides the discriminant of $\Chi $. 
We also know that the reduction of $\Chi $ modulo $(p_5')$ is
orthogonally stable. 
So the discriminant of $\Chi $ is in 
$$\{ 17 , 17 u, 17 p_7, 17 u p_7, 17 p_3p_5,  17 u p_3p_5, 17 p_7 p_3 p_5 , 
17 u p_7 p_3 p_5 \} .$$ 
Only for $d=17$ and $d=17up_7p_3p_5$ is the prime 2 not ramified 
in $K[\sqrt{d}]/K$. 
Now $p_5'$ is inert in this extension for $d=17$ and decomposed for
$d=17up_7p_3p_5= 68\sqrt{21} + 357$. As the reduction modulo $p_5'$ of $\Chi $ 
has orthogonal discriminant $O+$, we conclude that $d=68\sqrt{21} + 357$ is the only possibility. 

For the last character, number 33, we use a completely different argument. 
This character extends in two ways to the group $He\!\!:\!\!2$.
Take one of these extensions and restrict it to the maximal subgroup
$S_4(4)\!\!:\!\!4 \leq He\!\!:\!\!2$. 
This restriction is orthogonally stable and contains the following 
orthogonal simple characters: 
$$\begin{array}{|c|c|c|c|c|c|c|c|c|c|c|c|c|c|} 
	\hline
deg & 18ab & 68 & 50ab & 204 &170 & 153ab &  816 &  900 & 
	1020 & 256ab & 256  &  680 &  680  \\
	\hline
	K & \Q[i] & \Q & \Q[\sqrt{-2} ] & \Q& \Q & \Q[i] & \Q & 
	\Q & \Q & \Q[i] & \Q & \Q & \Q \\
	\hline
	mult & 1 & 1 & 1 & 2 & 1 & 1 & 4 & 7 & 6 & 3 & 1 & 2 & 5 \\
	\hline
	disc & 1 & 1 & 1 & . & -1 & -1 & . & 17 & . & 1 & 17 & . & 1 \\
	\hline
\end{array} 
$$
This table gives the degree of the orthogonal simple characters 
occurring in the restriction of the extension of $\Chi  $ to $S_4(4)\!\!:\!\!4$. 
The small letters indicate pairs of complex conjugate characters, 
whose character fields are  given in the next row. 
Their discriminant can hence be obtained by Proposition \ref{unitary}. 
The character of degree 68  is induced up from a
rational character of degree 34 of the group $S_4(4)\!\!:\!\!2$. 
Hence the invariant forms are orthogonal sums of two isometric forms 
and so the discriminant is a square. 
The same argument applies to the characters of degrees 170, 204, 816, 1020,
and the two characters of degree 680. 
For the character of degree 900 we apply \cite[Theorem 4.6]{DetChar} 
to conclude that the two characters of degree 450 of the group 
$S_4(4)\!\!:\!\!2$ have discriminant $-(17\pm 4\sqrt{17}) $ respectively. 
So the discriminant of the character of degree 900 is their product, 17. 
The characters of degree 256 restrict irreducibly to the simple group 
$S_4(4)$. The restriction is irreducible modulo 2 and 3, not orthogonally
stable modulo 17. So its discriminant is either 17 or 85. 
As the orthogonal discriminant of the reduction mod 2 is $O+$, 
we conclude that the discriminant is $\equiv 1 \pmod{8} $ (see Corollary \ref{mod2}) and hence it is 17.

\end{document}